\documentclass[a4paper,twoside,11pt]{article}

\usepackage{amsmath}
\usepackage{amsfonts}
\usepackage{amssymb}
\usepackage{amsthm}
\usepackage{amscd}

\input xy
\xyoption{all} \tolerance=500

\font\black=cmbx10 \font\sblack=cmbx7 \font\ssblack=cmbx5
\font\blackital=cmmib10  \skewchar\blackital='177
\font\sblackital=cmmib7 \skewchar\sblackital='177
\font\ssblackital=cmmib5 \skewchar\ssblackital='177
\font\sanss=cmss10 \font\ssanss=cmss8 scaled 900
\font\sssanss=cmss8 scaled 600 \font\blackboard=msbm10
\font\sblackboard=msbm7 \font\ssblackboard=msbm5
\font\caligr=eusm10 \font\scaligr=eusm7 \font\sscaligr=eusm5
 \font\fraktur=eufm10 \font\sfraktur=eufm7
\font\ssfraktur=eufm5

\font\bsymb=cmsy10 scaled\magstep2
\def\all#1{\setbox0=\hbox{\lower1.5pt\hbox{\bsymb
       \char"38}}\setbox1=\hbox{$_{#1}$} \box0\lower2pt\box1\;}
\def\exi#1{\setbox0=\hbox{\lower1.5pt\hbox{\bsymb \char"39}}
       \setbox1=\hbox{$_{#1}$} \box0\lower2pt\box1\;}

\def\tx#1{{\fam0\relax#1}}

\newfam\bifam
\textfont\bifam=\blackital \scriptfont\bifam=\sblackital
\scriptscriptfont\bifam=\ssblackital

\newfam\blfam
\textfont\blfam=\black \scriptfont\blfam=\sblack
\scriptscriptfont\blfam=\ssblack

\newfam\bbfam
\textfont\bbfam=\blackboard \scriptfont\bbfam=\sblackboard
\scriptscriptfont\bbfam=\ssblackboard

\newfam\ssfam
\textfont\ssfam=\sanss \scriptfont\ssfam=\ssanss
\scriptscriptfont\ssfam=\sssanss
\def\sss#1{{\fam\ssfam\relax#1}}

\newfam\clfam
\textfont\clfam=\caligr \scriptfont\clfam=\scaligr
\scriptscriptfont\clfam=\sscaligr

\newfam\frfam
\textfont\frfam=\fraktur \scriptfont\frfam=\sfraktur
\scriptscriptfont\frfam=\ssfraktur

\def\hpb#1{\setbox0=\hbox{${#1}$}
    \copy0 \kern-\wd0 \kern.2pt \box0}
\def\vpb#1{\setbox0=\hbox{${#1}$}
    \copy0 \kern-\wd0 \raise.08pt \box0}

\def\pmb#1{\setbox0\hbox{${#1}$} \copy0 \kern-\wd0 \kern.2pt \box0}
\def\pmbb#1{\setbox0\hbox{${#1}$} \copy0 \kern-\wd0
      \kern.2pt \copy0 \kern-\wd0 \kern.2pt \box0}
\def\pmbbb#1{\setbox0\hbox{${#1}$} \copy0 \kern-\wd0
      \kern.2pt \copy0 \kern-\wd0 \kern.2pt
    \copy0 \kern-\wd0 \kern.2pt \box0}
\def\pmxb#1{\setbox0\hbox{${#1}$} \copy0 \kern-\wd0
      \kern.2pt \copy0 \kern-\wd0 \kern.2pt
      \copy0 \kern-\wd0 \kern.2pt \copy0 \kern-\wd0 \kern.2pt \box0}
\def\pmxbb#1{\setbox0\hbox{${#1}$} \copy0 \kern-\wd0 \kern.2pt
      \copy0 \kern-\wd0 \kern.2pt
      \copy0 \kern-\wd0 \kern.2pt \copy0 \kern-\wd0 \kern.2pt
      \copy0 \kern-\wd0 \kern.2pt \box0}

\font\tengoth=eufm10 scaled\magstep1 \font\sevengoth=eufm7
\font\fivegoth=eufm5
\newfam\gothfam
\textfont\gothfam=\tengoth\scriptfont\gothfam=\sevengoth
  \scriptscriptfont\gothfam=\fivegoth

\mathchardef\za="710B  
\mathchardef\zb="710C  
\mathchardef\zg="710D  
\mathchardef\zd="710E  
\mathchardef\zve="710F 
\mathchardef\zz="7110  
\mathchardef\zh="7111  
\mathchardef\zvy="7112 
\mathchardef\zi="7113  
\mathchardef\zk="7114  
\mathchardef\zl="7115  
\mathchardef\zm="7116  
\mathchardef\zn="7117  
\mathchardef\zx="7118  
\mathchardef\zp="7119  
\mathchardef\zr="711A  
\mathchardef\zs="711B  
\mathchardef\zt="711C  
\mathchardef\zu="711D  
\mathchardef\zvf="711E 
\mathchardef\zq="711F  
\mathchardef\zc="7120  
\mathchardef\zw="7121  
\mathchardef\ze="7122  
\mathchardef\zy="7123  
\mathchardef\zf="7124  
\mathchardef\zvr="7125 
\mathchardef\zvs="7126 
\mathchardef\zf="7127  
\mathchardef\zG="7000  
\mathchardef\zD="7001  
\mathchardef\zY="7002  
\mathchardef\zL="7003  
\mathchardef\zX="7004  
\mathchardef\zP="7005  
\mathchardef\zS="7006  
\mathchardef\zU="7007  
\mathchardef\zF="7008  
\mathchardef\zW="700A  

\newcommand{\be}{\begin{equation}}
\newcommand{\ee}{\end{equation}}
\newcommand{\ra}{\rightarrow}

\newcommand{\bea}{\begin{eqnarray}}
\newcommand{\eea}{\end{eqnarray}}
\newcommand{\beas}{\begin{eqnarray*}}
\newcommand{\eeas}{\end{eqnarray*}}
\def\*{{\textstyle *}}
\newcommand{\R}{{\mathbb R}}

\newcommand{\Z}{{\mathbb Z}}
\newcommand{\N}{{\mathbb N}}

\newcommand{\we}{\wedge}
\newcommand{\nn}{\nonumber}
\newcommand{\ot}{\otimes}

\newcommand{\pa}{\partial}
\newcommand{\ti}{\times}
\newcommand{\A}{{\cal A}}

\newcommand{\ad}{{\rm ad}}
\newcommand{\rS}{]^{SN}}

\newcommand{\X}{{\cal X}}

\newcommand{\Ll}{{\pounds}}
\def\lna{\lbrack\! \lbrack}
\def\rna{\rbrack\! \rbrack}

\def\lan{\langle}
\def\ran{\rangle}

\def\ix{\operatorname{i}}
\def\Ker{\operatorname{Ker}}
\def\im{\operatorname{Im}}

\def\V{{\cal V}}

\def\cD{{\cal D}}

\def\cE{{\cal E}}
\def\cF{{\cal F}}

\def\cL{{\cal L}}

\def\cM{{\cal M}}
\def\cP{{\cal P}}
\def\cR{{\cal R}}

\def\cT{{\cal T}}

\def\wh{\widehat}
\def\wt{\widetilde}

\def\Sec{\sss{Sec}}

\def\la{\langle}
\def\ran{\rangle}

\newcommand{\fG}{\mathfrak{G}}

\def\sT{{\sss T}}

\def\st{{\sss t}}

\def\xd{\tx{d}}
\def\xi{\tx{i}}

\def\k{{\mathfrak{k}}}
\def\End{\operatorname{\sss End}}
\def\Der{\operatorname{\sss Der}}
\def\Diff{\operatorname{\sss Diff}}
\def\Gr{\mathbb{A}}
\def\Ts{\operatorname{\sss Ts}}
\def\Al{\operatorname{\sss Al}}
\def\im{\operatorname{Im}}
\def\gl{{\mathfrak{gl}}}
\def\g{{\mathfrak g}}
\def\1{\mathbf{1}}

\newcommand{\rF}{]^{FN}}

\def\H{{\cal G}}
\def\V{{\cal V}}
\def\J{{\cal J}}

\newcommand{\bk}[2]{\ensuremath{\langle #1 | #2 \rangle}}
\newcommand{\pr}[2]{\ensuremath{\langle #1 , #2 \rangle}}
\newcommand{\fA}{\mathfrak{A}}
\newcommand{\fH}{\mathfrak{H}}
\newcommand{\half}{{\frac{1}{2}}}
\newcommand{\third}{{\frac{1}{3}}}

\def\on{\operatorname}
\def\xdp{{\on{d}^\Pi}}
\newcommand{\dd}{\mathrm{d}}
\def\Id{\on{Id}}

\def\M{{\sss M}}

\def\y{\mathbf{y}}
\def\QD{\operatorname{\sss QD}}

\newtheorem{thm}{Theorem}[section]
\newtheorem{theorem}{Theorem}[section]
\newtheorem{proposition}[thm]{Proposition}

\newtheorem{corollary}[thm]{Corollary}

\theoremstyle{definition}
\newtheorem{example}[thm]{Example}
\newtheorem{remark}[thm]{Remark}
\newtheorem{definition}[thm]{Definition}

\begin{document}
\pagestyle{myheadings}
\markboth{Janusz Grabowski} {Brackets}

%
%

\title{BRACKETS\thanks{Mini-course held at XXI Fall Workshops on Geometry and Physics, Burgos (Spain), 2012.}}

\author{Janusz Grabowski\thanks{Research
supported by the Polish Ministry of Science and Higher Education under the grant N N201 416839.}\\ \\
Institute of Mathematics, Polish Academy of Sciences\\ \'Sniadeckich 8, 00-956 Warszawa, POLAND\\ {\tt jagrab@impan.pl}}

\date{}
\maketitle
\begin{center}
\emph{Dedicated to the memory of Jean-Louis Loday}
\end{center}

\begin{abstract}
 We review origins and main properties of the most important bracket operations appearing canonically in differential geometry and mathematical physics in the classical, as well as in the supergeometric setting. The review is supplemented by some new concepts and examples.
\end{abstract}

\section{Introduction}
In algebra, `brackets' are usually understood as non-associative
operations on vector spaces or modules, with Lie brackets as the
main example. The aim of this notes is to present a survey
of bracket operations playing an important role in geometry and
physics applications. What we will consider are mainly canonical Lie and, more
generally, Loday brackets in the standard, as well as
superalgebraic setting. Among them are Poisson and Jacobi (more
generally, Kirillov) brackets,  Schouten-Nijenhuis,
Nijehuis-Richardson, and Fr\"olicher-Nijenhuis brackets, Lie
algebroid brackets, Courant and Dorfman brackets, $n$-ary brackets of Filippov and Nambu,
etc.

A proper understanding of the roots and properties of all these
brackets requires a basic knowledge of superalgebra and graded
differential geometry whose rudiments will be also outlined in these notes. We will add also a few new concepts and examples. Of
course, the subject is so extensive that we are only able to sketch
selected problems and cite only a small part of the existing literature.
We hope, however, that this review could be of some interest for those who encounter the brackets in their work with problems of contemporary mathematics and physics.

\section{Lie and graded Lie algebras}

\subsection{Algebras} By an \emph{algebra} on a vector space $\A$ over a field
$\k$ we will understand a bilinear operation on $\A$, \be\A\ti
\A\ni(x,y)\mapsto x\circ y\in \A\,. \ee
In most cases we can consider as well algebras over commutative rings.
The algebra $(\A,\circ)$ we call
\emph{commutative} if the operation is commutative, $x\circ
y=y\circ x$, and \emph{anti-commutative} or \emph{skew-symmetric}
if $x\circ y=y\circ x$. We call $(\A,\circ)$ \emph{unital} if it
has a \emph{unit}, i.e.  an element $\1\in\A$ for which $\1\circ
x=x\circ\1=x$. Note that the unit is unique if it exists. A
commutative associative operation we will usually denote
"$\cdot$", or even write it simply as the juxtaposition, e.g. 
$xy$.

If $\{ x_i\}_{i\in I}$ is a basis of $\A$, the algebra
structure is uniquely determined by the \emph{structure constants}
$c^k_{ij}$, where (the summation convention is used)
\be x_i\circ x_j=c^k_{ij}x_k\,. \ee

\begin{example} Let us observe that any vector space $V$ gives rise to a
canonical nontrivial algebra structure. Namely, the space
$\A=\gl(V)=\End_\k(V)$ of all linear endomorphisms $x:V\to V$ is
an algebra with the operation being just the composition of maps.
This algebra is \emph{associative}, i.e.  the map \be
m:\A\to\gl(\A)\,,\quad m_x(y)=x\circ y\,,\ee is an algebra
homomorphism,
\begin{equation}\label{1}m_{x\circ y}=m_x\circ m_y\,.\end{equation} Here, of course, the first "$\circ$" is the operation in $\A$ and the second in $\gl(\A)$. In other words, the product in
$\A$ satisfies the identity
\begin{equation}\label{ass} (x\circ y)\circ z=x\circ(y\circ z)\,.\end{equation} We call $m$ the \emph{(left) regular representation} of $(\A,\circ)$.
\end{example}

\begin{example} If $M$ is a topological space, then the set $C(M)$ of all real continuous
functions on $M$ is canonically a commutative associative algebra
over $\R$ with the point-wise multiplication. Similarly, if $M$ is
a smooth manifold, then the set $C^\infty(M)$ of all real smooth
functions on $M$ is also a commutative associative algebra.
\end{example}

Having one binary operation, we can easily produce other
operations. For instance, we can consider the \emph{commutator}
$[x,y]=x\circ y-y\circ x$, or the \emph{anti-commutator}
(\emph{symmetrizer}) $x\vee y=x\circ y+y\circ x$. For linear
operators, this produces canonical, this time a skew-symmetric
(resp., symmetric), operation in $\A=\gl(V)$ which is in general no
longer associative; for the symmetric product we have only a week
version of the associativity:
\begin{equation}\label{jo}
(x\vee y)\vee(x\vee x)=x\vee(y\vee(x\vee x))\,.
\end{equation} These structures are prototypes of what we call a \emph{Lie algebra} or, respectively, a \emph{Jordan algebra}. For instance,
we can easily check the following analog of (\ref{1}):
\begin{equation}\label{2} \ad_{[x,y]}=[\ad_x,\ad_y]\,,\end{equation} where $\ad:\A\to\gl(\A)$ is the corresponding regular representation for the commutator, $\ad_x(y)=[x,y]$. In
other words, $\ad$ is a homomorphism of the brackets, i.e. 
\begin{equation}\label{J} [[x,y],z]=[x,[y,z]]-[y,[x,z]]\,.\end{equation} Identity (\ref{J}) we call the \emph{Jacobi identity}.

\begin{remark} Let us remark that sometimes by the Jacobi identity one understands the identity
\begin{equation}\label{J1} [x,[y,z]]+[y,[z,x]]+[z,[x,y]]=0\,.\end{equation} For a skew-symmetric operation (bracket) both versions are equivalent, but for brackets which are not
skew-symmetric this is no longer true. The advantage of (\ref{J}) is its clear algebraic meaning:  $\ad$ is a representation; therefore the Jacobi identity will
be for us always (\ref{J}).
\end{remark}

The Jacobi identity means, equivalently, that operators $\ad_x$ are
derivations of the bracket. Recall that a \emph{derivation} of an
algebra $(\A,\circ)$ is a map $D\in\gl(\A)$ such that, for all
$x,y\in\A$,
\begin{equation}\label{der}
D(x\circ y)=D(x)\circ y+x\circ D(y)\,,
\end{equation}
i.e.  the \emph{Leibniz rule} is satisfied.
A trivial but very useful observation giving a method of constructing derivations is that, if $\A$ is freely generated by $(x_i)$ and $y_i\in\A$, then there
is a unique derivation $D$ of $\A$ such that $D(x_i)=y_i$.

A bracket $[\cdot,\cdot]$ satisfying the Jacobi identity is called a
\emph{Leibniz bracket} or \emph{Loday bracket} and the
corresponding algebra a \emph{Leibniz (Loday) algebra}.
This terminology goes back to the work of J.-L.~Loday who discovered that one can skip
the skew-symmetry assumption in the definition of a Lie algebra,
still having a possibility to define appropriate (co)homology (see
\cite{Lo1,LP} and \cite[Chapter 10.6]{Lo}). Loday himself called these structures \emph{Leibniz algebras}.

Of course, if the Loday bracket is additionally skew-symmetric, we speak about a \emph{Lie
bracket} and a \emph{Lie algebra}. The space $\A=\gl(V)$ of linear
operators with the commutator bracket is therefore a canonical
example of a Lie algebra.

Similarly, a vector space equipped with a symmetric operation
$\vee$ satisfying (\ref{jo}) we call a \emph{Jordan algebra}. The
space $\A=\gl(V)$ of linear operators with the anti-commutator
bracket is therefore a canonical example of a Jordan algebra.
Operators symmetric with respect to a certain anti-involution
$x\mapsto x^\dag$ (think of Hermitian operators) form a Jordan
subalgebra in $\gl(V)$.

If $V$ has additionally an algebra structure with respect to a
product "$\circ$", we can distinguish canonically a Lie subalgebra
in $\A=\gl(V)$. This is one of the major ways of obtaining Lie algebra
structures.

\begin{proposition} For any algebra $(\A,\circ)$, the space $\Der(\A,\circ)$ of derivations of $\A$ is a Lie subalgebra in $\gl(\A)$ with respect to the commutator
bracket.
\end{proposition}

\noindent If the product "$\circ$" is fixed, instead of
$\Der(V,\circ)$ we will write simply $\Der(V)$.

\begin{example} Let $M$ be a manifold and let $\A=C^\infty(M)$ be the commutative associative algebra of
smooth functions on $M$. Then, the Lie algebra $\Der(\A)$ is canonically
identified with the Lie algebra $\X^1(M)$ of smooth vector fields
on $M$ with the Lie bracket of vector fields.
\end{example}

\noindent{\bf Problem.} Show that there are no non-zero
derivations of the algebra $C(\R)$ of all continuous functions on
$\R$, i.e.  the differential calculus for $\A=C(\R)$ is trivial.


\subsection{Modules}
Having a (commutative) associative, (resp., Leibniz, Lie, Jordan,
etc.) algebra $(A,\circ)$ over $\k$, we define its \emph{module}
as a vector space $V$ over $\k$ equipped with two operations, $A\ti
V\ni(a,v)\mapsto a\circ v\in V$ and $V\ti A\ni(v,a)\mapsto v\circ
a\in V$, such that $A\oplus V$ becomes also a (commutative)
associative, (resp., Leibniz, Lie, Jordan, etc.) algebra with
obviously defined operation, denoted with some abuse of notation
also $\circ$, which is trivial on $V$, $v_1\circ v_2=0$.

\begin{example} If $(A,\cdot)$ is a commutative algebra, then its module $V$ is defined by a multiplication
$A\ti V\ni(a,v)\mapsto av\in V$ such that $(a_1a_2)v=a_1(a_2v)$.
The other multiplication $V\ti A\ni(v,a)\mapsto va\in V$ is
uniquely determined by the symmetry, $va=av$.
\end{example}

\begin{example}
If $\zt:E\to M$ is a vector bundle, then the space $\cE=\Sec(E)$
of all sections of $E$ is canonically an $\A=C^\infty(M)$-module
with the obvious multiplication ($\A$ is commutative, so the left
and the right multiplications are equal) of section by
functions on $M$.
\end{example}

\begin{example}
In particular, the space $\X^1(M)$ of vector fields on $M$ is a
$C^\infty(M)$-module and the canonical Lie bracket of vector
fields is related to this module structure by the following
`Leibniz rule':
\begin{equation}\label{anc} [X,fY]=f[X,Y]+X(f)Y\,.
\end{equation}
Actually, the action of vector fields on functions makes $C^\infty(M)$ into a $\X^1(M)$-module, so that $X(f)$ can be viewed as a bracket, $[X,f]$, which makes the space \be
X^1(M)\oplus
C^\infty(M)=\Sec(\sT M\ti\R)\ee
of linear first-order differential operators on $M$ into a Lie
algebra.

Note that the isomorphism class of the Lie algebra $\X^1(M)$ completely determines $M$ up to a diffeomorphism, exactly like does it the associative algebra $C^\infty(M)$ \cite{Gra1,Gra2,SP}. Similar results are valid also for the Lie algebras of first-order differential operators \cite{GP}, Kirillov's local Lie algebras \cite{Gra4}, and for supermanifolds \cite{GKoP}.
\end{example}

\subsection{Graded algebras}
Let $K$ be a commutative associative ring with identity, $U(K)$ be the group of invertible elements
of $K$, and let $G$ be a commutative semigroup. A map $\ze: G \ti G \to U(K)$ is called a \emph{factor} on $G$ if
\be\ze(g,h)\ze(h,g) = 1\,,\quad p(g)=\ze(g,g) = \pm 1\,,\quad \text{and}\quad\ze(f,g+h)=\ze(f,g)\ze(f,h)\,,\ee
for all $f, g, h\in G$.
Let $V$ be a $G$-graded $K$-algebra, $V=\oplus_{g\in G}V^g$. Elements $x$ from $V^g$ we call \emph{homogeneous
of degree}  $g$ (or \emph{of weight $g$}) and denote $g=w(x)$. The algebra $V$ is called
\emph{$\ze$-commutative} if
\be a\circ b = \ze(w(a), w(b))b\circ a\ee
for all $G$-homogeneous elements $a, b \in V$. Homogeneous elements $a$ with $p(w(a))=-1$
we call \emph{odd}, the other homogeneous elements we call \emph{even}.
In what follows, $K$ will be $\R$ and $\ze$ will take the form $\ze(g,h)=(-1)^{\la g| h\ran}$ for
a `scalar product' $\la\cdot |\cdot\ran:G\ti G\to\Z$, with
$G=\Z_2$ (superalgebra case), $G=\Z^n$, or $G=\N^n$.
This means that we use the factor as a \emph{sign rule} which can be applied (separately) to
any axiom: we change the sign by $(-1)^{\la g|h\ran}$ whenever
two consecutive homogeneous elements $x\in V^g$ and $y\in V^h$ are
interchanged.

We say that our operation "$\circ$" is \emph{of degree $k$} (or simply
\emph{even} or \emph{odd} for $k=0,1\in G=\Z_2$) if
\begin{equation}\label{deg}V^g\circ V^h\subset V^{g+h+k}\,.\end{equation}
For operations of degree $k$ it is natural to consider the sin rules in the form
\be a\circ b=(-1)^{\la w(a)+k|w(b)+k\ran} b\circ a\,.\ee
If $G$ is a group, we can always work with degree 0, making the
corresponding shift, $V[k]$, in the grading, where
$V[k]^i=V^{i-k}$.
For $\Z^n$-gradings, we will use by default $\la i|j\ran=i_1j_1+\cdots  i_nj_n$, but other bi-additive pairings
$\Z^n\ti\Z^n\to\Z$ are also acceptable. In this way we get the
graded versions (of degree $k\in\Z^n$) of our structures. For
instance, the graded symmetry of degree $k$ reads
\begin{equation}\
\label{s1}x\circ y=(-1)^{\la w(x)+k|w(y)+k\ran}y\circ x
\,,\end{equation} the graded skew-symmetry of degree $k$ reads
\begin{equation}\
\label{s2}x\circ y=-(-1)^{\la w(x)+k|w(y)+k\ran}y\circ x
\,,\end{equation} and the graded Jacobi identity of degree $k$
\begin{equation}\label{gJ}[[ x,y],z]=[ x,[ y,z]]-(-1)^{\lan w(x)+k|w(y)+k\ran}[ y,[ x,z]]\,.\end{equation} If $k=0$, we speak simply about the graded commutativity, the graded skew-symmetry, the graded Jacobi identity,
etc. The above identities for the degree $k$ are just the graded
commutativity, graded skew-symmetry and the graded Jacobi identity
for the shifted grading. Note that no sign appears in the
associativity property (\ref{ass}).

\begin{example}
Starting with a vector bundle $\zt:E\to M$ we can consider the
\emph{tensor algebra} of its sections, \be\Ts(E)=\Sec(\ot
E)=\oplus_{i=0}^\infty\Sec(E^{\ot i})\,,\ee which is clearly an
associative (but noncommutative) graded algebra with respect to
the tensor product. It contains the \emph{Grassmann algebra}
\be\Gr(E)=\oplus_{i=0}^\infty\Gr^i(E)=\oplus_{i=0}^\infty\Sec(\we^i
E)\ee which consists of skew-symmetric tensors and is a
graded-commutative associative algebra with respect to the wedge
product $\wedge$, the skew-symmetrization of the tensor product.
In particular, for $E=\sT M$ and $E=\sT^\*M$, we obtain the $\N$-graded commutative associative
algebras $\X(M)=\oplus_{i=0}^\infty\X^i(M)$ and
$\Omega(M)=\oplus_{i=0}^\infty\zW^i(M)$ of multivector fields and
differential forms, respectively.
\end{example}

Given a graded algebra $(A=\oplus_{i\in{\Z^n}}A^i,\circ)$, we
define the space $\Der^k(A)$ of \emph{graded derivations} of
degree $k\in\Z^n$ as the space of maps $D:A\to A$ of degree $k$ such that
\begin{equation}\label{gDer} D(x\circ y)=D(x)\circ y+(-1)^{\la w(x)|k\ran}x\circ D(y)
\end{equation} is satisfied. Then, the graded space of (graded) derivations $\Der(A)=\oplus_{k\in \Z^n}\Der^k(A)$ is a graded
Lie algebra with respect to the graded commutator,
\begin{equation}\label{gcom} [D_1,D_2]=D_1D_2-(-1)^{\la k_1|k_2\ran}D_2D_1\,,
\end{equation} where $D_i\in\Der^{k_i}(A)$, $i=1,2$.

\medskip\noindent{\bf Problem.\ } What are graded derivations of degree 0 of the Grassmann algebra $\zW(M)$ of differential forms?

\medskip
Note that for a skew-symmetric bracket operation
$[\cdot,\cdot]:V\we V\to V$ on a, say, finite-dimensional vector
space $V$, its dual is a certain map $\xd:V^\*\to V^\*\we V^\*$.
As the Grassmann algebra $\Gr(V^\*)$ is freely generated by any
basis of $V^\*$, this maps gives rise to a uniquely defined graded
derivation $\xd:\Gr(V^\*)\to\Gr(V^\*)$ of degree 1. Indeed, we can
inductively define \be\xd(v_0\we\cdots\we v_n)=\xd v_0\we
v_1\we\cdots\we v_n-v_0\we\xd (v_1\we\cdots\we v_n)\,.\ee
Explicitly,
\begin{equation}\label{CE}\xd\za(x_0,\dots,x_n)=\sum_{i<j}(-1)^{i+j}\za([x_i,x_j],x_0,\dots,\hat{x_i},\dots,\hat{x_j}\dots,x_n)\,.
\end{equation}
If $V$ is a Lie algebra, this derivation is a \emph{homological} operator, $\xd^2=0$,
called the \emph{Chevalley-Eilenber cohomology operator} and defining the \emph{Lie algebra
cohomology} in the standard way:
$H^\bullet(V)=(\Ker \xd/\im \xd)^\bullet$. The cohomology operator $\xd$ contains the full
information about the Lie algebra structure, as the bracket is its
dual $\xd^\*$. We can therefore formulate
an equivalent definition of a finite-dimensional Lie algebra as
follows.

\begin{proposition} A Lie algebra structure on a
finite-dimensional vector space $V$ is a degree 1 derivation of the Grassmann algebra $\Gr(V^\*)$
which is homological, i.e.  $d^2=0$.
\end{proposition}

We will formulate later a similar fact for Lie algebroids
understanding, after Vaintrob \cite{Va,Va1}, the corresponding derivations as \emph{homological vector fields}.

\subsection{Gerstenhaber and Nijenhuis-Richardson brackets}

We denote with $\M^p(V)$ the space of all $p$-linear maps $A :
V^{p}\to V$ if $p> 0$. We put $\M^{0}( V) = V$ and we set
$\M(V)= \oplus_{p\ge 0} \M^p( V)$.
On the graded vector space $\M( V)$ we define the operation $i: \M(
V)^2\to\M( V)$ of degree $-1$ by: $i(B)A=0$ if $A\in\M^{0}(V)$, and
\bea\nn&i(B)A(x_1,\dots,x_{ a+b-1}) =\\
&\sum_{k=1}^a(-1)^{(k-1)(b-1)}A(x_1,\dots,x_{k-1},B(x_{k}\dots,x_{k+b-1}),x_{k+b},\dots,x_{a+b-1})\label{Ge}
\eea if $A\in\M^{a}( V)$, $a>0$, and $B\in\M^b( V)$.
Define now the
bracket $[\cdot,\cdot]^G: \M( V)^2\to\M( V)$ of degree -1 by
\begin{equation}\label{Ge1} [A,B]^G=i(B)A- (-1)^{(a-1)(b-1)} i(A)B\,,\quad A\in\M^a( V)\,,\quad B\in\M^b( V)\,.
\end{equation} This bracket is an extension of the usual commutator bracket in $M^1(V)=\gl(V)$, called the \emph{Gerstenhaber bracket} \cite{Ge}.

For the graded subspace  $\Al( V)$ of $\M(V)$ of alternating (skew-symmetric)
mappings, define the bracket $[\cdot,\cdot]^{RN} : \Al( V)^2\to \Al(
V)$ of degree -1 by
\begin{equation}\label{aNR} [A,B]^{NR}=\frac{(a+b- 1)!}{a\,!\,b\,!}\zs([A,B]^G)\,,\quad
A\in\Al^a( V)\,,\quad B\in\Al^b( V)\,,
\end{equation} where $\zs$ stands for the anti-symmetrization projector in $\M( V)$. This bracket is called the \emph{(algebraic) Nijenhuis-Richardson bracket}. The importance of the above brackets indicates the following observation which shows that they serve for determining associative and Lie algebra structures, together with the corresponding cohomology.

\begin{proposition} The brackets $[\cdot,\cdot]^G $ and $[\cdot,\cdot]^{NR}$ are graded Lie brackets of degree $-1$ on $\M( V)$ and $\Al( V)$, respectively.
Moreover, a map $A\in\M^2(V)$ (resp., $A\in\Al^2(V)$) defines an
associative (resp., Lie) algebra structure on $V$ if and only if
$A$ is a homological element, i.e.  $[A,A]^G=0$ (resp., $[A,A]^{NR}=0$). In
this case, the adjoint map $\pa_A:\M(V)\to\M(V)\,, \pa_A(B)=[B,A]^G$ (resp.,
$\pa_A:\Al(V)\to\Al(V)\,, \pa_A(B)=[B,A]^{NR}$) is homogeneous of
degree 1 and satisfies $\pa_A^2=0$, so that it defines a cohomology, called the \emph{Hochschild} (resp., \emph{Chevalley-Eilenberg}) \emph{cohomology}.
\end{proposition}

\subsection{Poisson brackets}
If we have an isomorphism of the vector bundles $\sT^\*M$ and $\sT
M$, thus inducing an isomorphism
$\zW^1(M)\ni\za\mapsto\wh{\za}\in\X(M)$, then we can transform the
Lie bracket of vector fields into the space $\zW^1(M)$ of
one-forms such that
\begin{equation}\label{lbt}[\za,\zb]^{\,\wh{}}=[\wh{\za},\wh{\zb}]\,.
\end{equation} For instance, a symplectic form $\zw$ on $M$ induces such an isomorphism, $\wt{\zw}:\sT M\to\sT^\* M$, and the corresponding bracket
$[\za,\zb]_\zw$ \emph{via} \be\zw(\cdot,\wh{\za})=\za\,.\ee As
easy calculations show, for any one-forms $\za,\zb$ and any vector
field $X$, \be
0=\xd\zw(\wh{\za},\wh{\zb},X)=-\xd\za(\wh{\zb},X)+\xd\zb(\wh{\za},X)-\zw([\wh{\zb},\wh{\za}],X)+
i_X\xd(\zw(\wh{\za},\wh{\zb}))\,,\ee so that
\begin{eqnarray}\nn[\za,\zb]_\zw&=&i_{\wh{\za}}\xd\zb-i_{\wh{\zb}}\xd\za+\xd(\zw(\wh{\za},\wh{\zb}))\\
&=&\Ll_{\wh{\za}}\zb-\Ll_{\wh{\zb}}\za-\xd(\zw(\wh{\za},\wh{\zb}))\,,\label{Kb}
\end{eqnarray} where $\Ll$ denotes the Lie derivative. The bracket (\ref{Kb}) is called the \emph{Koszul bracket} of
one-forms. If $\za$ and $\zb$ are exact, $\za=\xd f$ and $\zb=\xd
g$, the vector fields $X_f=\wh{\xd f}$ and $X_g=\wh{\xd g}$ are
called the \emph{Hamiltonian vector fields} with \emph{Hamiltonians}
$f$ and $g$, respectively. In this case, we have \be[\xd f,\xd
g]_\zw=\xd(\zw(X_f,X_g))=\xd (X_f(g))\ee and \be[X_f,X_g]=[\wh{\xd
f},\wh{\xd g}]=\wh{[\xd f,\xd
g]}={\xd\left(\zw(X_f,X_g)\right)}^{\widehat{\ \ }}\,,\ee so that the de Rham
derivative is a homomorphism of the bracket
\begin{equation}\label{Pb}\{ f,g\}_\zw=\zw(X_f,X_g)=X_f(g)
\end{equation} on $C^\infty(M)$ into the Koszul bracket,
\be\xd\{ f,g\}_\zw=[\xd f,\xd g]_\zw\,.\ee Actually, $\{
f,g\}_\zw$ is a Lie bracket:
\bea\nn\{\{ f,g\}_\zw,h\}_\zw&=&X_{\{ f,g\}_\zw}(h)=[X_f,X_g](h)=X_f(X_g(h))-X_g(X_f(h))\\
&=&\{ f,\{ g,h\}_\zw\}_\zw-\{ g,\{ f,h\}_\zw\}_\zw\,. \eea

\begin{definition}
A Lie bracket $[\cdot,\cdot]$ on an associative algebra
$(V,\circ)$ such that the operators $\ad_x$ act as derivations
also for the associative multiplication, i.e.  the Leibniz rule
\begin{equation}\label{lrpa} [x,y\circ z]=[x,y]\circ z+ y\circ[x,z]
\end{equation} is satisfied, is called a \emph{Poisson bracket}, and the triple $(V,\circ,[\cdot,\cdot])$ a \emph{Poisson
algebra}.
\end{definition}

Note that any associative algebra is automatically a
Poisson algebra with respect to the commutator bracket. Of course, this bracket is trivial for any commutative algebra,
so Poisson brackets are extra structures for the latter.

\begin{example} If $(M,\zw)$ is a symplectic manifold, then the bracket (\ref{Pb}) is a Poisson bracket and turns
$C^\infty(M)$ into a Poisson algebra. Indeed, \be\{
f,gh\}_\zw=X_f(gh)=X_f(g)h+gX_f(h)=\{ f,g\}_\zw h+g\{
f,h\}_\zw\,.\ee
\end{example}

Due to the Leibniz rule (\ref{lrpa}), any Poisson bracket
$\{\cdot,\cdot\}$ on $C^\infty(M)$ is represented by a bivector
field $\zL$, so that
\begin{equation}\label{pst}
\{ f,g\}=\{ f,g\}_\zL=\la\zL,\xd f\we\xd g\ran\,.
\end{equation} This is the contravariant version of (\ref{Pb}). Of course, in view of the Jacobi identity, the tensor $\zL$ must
satisfy an additional condition. In local
coordinates, if \be\zL =
\dfrac{1}{2}\zL^{ij}(x)\dfrac{\partial}{\partial x^i }\wedge
\dfrac{\partial}{\partial x^j }\,, \ee then \be \{
f,g\}_\zL=\zL^{ij}(x)\frac{\pa f}{\pa x^i}\frac{\pa g}{\pa x^j}\ee
is a Poisson bracket if and only if, for all $j,k,l$,
\begin{equation}\label{ptnsr}\sum_i(\zL^{ij}\dfrac{\partial \zL^{kl}}{\partial x^i } +
\zL^{ik}\dfrac{\partial \zL^{lj}}{\partial x^i }  +
\zL^{il}\dfrac{\partial\zL^{jk}}{\partial x^i }) =0 \,.
\end{equation}
Such tensors $\zL$ we call \emph{Poisson tensors} or
\emph{Poisson structures}. The above conditions have a nice interpretation in terms of the so called \emph{Schouten-Nijenhuis bracket} (see the
next paragraph).

\begin{remark}
It can be proven that the skew-symmetry of Poisson brackets on $C^\infty(M)$
follows from the Leibniz rule and the Jacobi identity
\cite{GMbin}, so it is a superfluous condition in the definition. In \cite{Gra}, a canonical extension of the Poisson bracket of
functions to a graded Lie bracket on differential forms has been
constructed. This bracket, however is not a graded Poisson bracket,
as the Leibniz rule is not satisfied. Actually, it is a
second-order bracket.
\end{remark}

\begin{example} (KKS-structure) Let $\g$ be a finite-dimensional
real Lie algebra with a Lie bracket $[\cdot,\cdot]$ and let
$c^k_{ij}$ be the structure constants with respect to a basis
$x_1,\dots,x_n$. Note that $x_1,\dots,x_n$ can be viewed as linear
functions defining a coordinate system on the dual space $\g^\*$.
Then, there is a uniquely determined Poisson bracket
$\{\cdot,\cdot\}$ on $\g^\*$ (\emph{Kostant-Kirillov-Souriau bracket}) such that \be\{
x_i,x_j\}=[x_i,x_j]=c^k_{ij}x_k\,.\ee Indeed, it is easy to see
that the corresponding tensor must be
$\zL=\frac{1}{2}c^k_{ij}x_k\pa_{x_i}\we\pa_{x_j}$ which satisfies
(\ref{ptnsr}), as the latter is in this case equivalent to the
Jacobi identity for $[\cdot,\cdot]$. The Poisson tensor is
\emph{linear} in the obvious sense and hence the corresponding Poisson
bracket is closed on polynomial functions.
\end{example}

One important observation is that the above correspondence between
Lie brackets on a vector space and linear Poisson tensors on the
dual remains valid for an arbitrary vector bundle $\zt:E\to M$
(see Theorem \ref{ty}). Of course, linear functions on a vector
bundle are understood as functions which are linear along fibres, and the
linearity of a Poisson tensor means that the corresponding Poisson
bracket is closed on linear functions. Automatically, it is closed
on the space of polynomial functions which becomes, in this way, a
Lie algebra.

\begin{example} It is well known that the cotangent bundle $\sT^\* M$ posses a canonical symplectic structure $\zw_M$,
thus $C^\infty(\sT^\* M)$ is canonically a Poisson algebra. There
are local affine coordinates $(q^a,p_a)$ on $\sT^\* M$, called
\emph{Darboux coordinates}, in which $\zw_M=\xd p_a\we\xd q^a$ and
in which this Poisson bracket reads as
\begin{equation}
\{   f,g\}   =   \frac{\partial    f}{\partial p_a} \frac{\partial
g}{\partial q^a}-\frac{\partial  f}{\partial q^a}\frac{\partial
g}{\partial p_a}\,.\label{b}
\end{equation}
The corresponding Poisson tensor $\zL_M=\pa_{p_a}\we\pa_{q^a}$ is
the `inverse' of the symplectic form $\zw_M$ in the sense that,
\emph{via} the contraction, it defines the inverse isomorphism
\begin{equation}\label{beta}\wt{\zL_M}:\sT^\*\sT^\* M\to\sT\sT^\* M\,.
\end{equation}
The tensor is linear, because the bracket is closed on linear
(i.e.  linear in $p$'s) functions. Consequently,
polynomial (in $p$'s) functions form a Lie algebra. This is the Lie algebra of symbols of  differential operators on $M$.
\end{example}

\begin{theorem}(Darboux Theorem)
Each symplectic Poisson bracket can be written locally in the
form (\ref{b}).
\end{theorem}

Let us observe that the linear Poisson bracket (\ref{b}) is
\emph{de facto} equivalent to the Lie bracket of vector fields on
$M$. Indeed, we can identify any vector field $X$ on $M$ with the
corresponding linear function $\zi_X$ on $\sT^\* M$ in an obvious
way: $\zi(X)(\za_q)=\la X(q),\za_q\ran$.  In local coordinates,
\begin{equation}\label{zi}\zi({f^a(q)\pa_{q^a}})=f^a(q)p_a\,.
\end{equation}
It is easy to see now that
\begin{eqnarray}\label{i1} \{\zi({X}),\zi(Y)\}&=&\zi({[X,Y]})\,,\\
\{\zi({X}),f\}&=&X(f)\,,\label{i2}
\end{eqnarray} where $f$ is any basic function on $\sT^\*M$ interpreted as a function on $M$.
A more detailed study of Poisson brackets and related structures can be found in \cite{Vs}.

\subsection{Jacobi and graded Jacobi brackets}
A construction of a Lie bracket, similar to that for
functions on a symplectic manifold, can be done in the case of a
contact manifold $(M,\za)$. We call this bracket the
\emph{Legendre bracket}.

\begin{example} Not going into a general theory, let us recall that any contact
form can be locally written as $\za=\xd z-p_a\xd q^a$ (Darboux
theorem) 
and the Legendre bracket of functions $f,g$ on $M$ in these
coordinates reads as
\begin{equation}\label{cobr}\{ f,g\}_\za =\frac{\pa f}{\pa{p_a}}\frac{\pa g}{\pa
q^a}-\frac{\pa f}{\pa q^a}\frac{\pa g}{\pa{p_a}}+\frac{\pa f}{\pa
z}\left(g-p_a\frac{\pa g}{\pa{p_a}}\right)-\left(f-\frac{\pa
f}{\pa{p_a}}p_a\right)\frac{\pa g}{\pa z}\,.
\end{equation}
This bracket is not Poisson, since the Leibniz rule is not
satisfied: the operators $\{ f,\cdot\}_\za$ act on $C^\infty(M)$
as first-order differential operators, not derivations. This can
be expressed in terms of a \emph{generalized Leibniz rule}:
\begin{equation}\label{Jb}
\{ f,gh\}_\za=\{ f,g\}_\za h+g\{f,h\}_\za-\{ f,\1\}_\za gh\,.
\end{equation}
\end{example}

A Lie bracket on a (commutative) associative unital algebra,
satisfying (\ref{Jb}), we call a \emph{Jacobi bracket}. Thus the
bracket (\ref{cobr}) is an example of a Jacobi bracket on
$C^\infty(M)$. In general, Jacobi brackets on $C^\infty(M)$ are represented by
pairs $(\zL,\zG)$, where $\zL$ is a bivector field and $\zG$ is a
vector field, by
\begin{equation}\label{Jabra}\{ f,g\}_{(\zL,\zG)}=\zL(\xd f,\xd g)+\zG(f)g-f\zG(g)\,.
\end{equation}
The pair $(\zL,\zG)$ is called a \emph{Jacobi structure} \cite{Li}. For the Legendre bracket (\ref{cobr}),
\be\zL=\pa_{p_a}\we\pa_{q^a}+p_a\pa_{p_a}\we\pa_z\,,\quad
\zG=\pa_z\,,\ee Of course, the Jacobi identity for the bracket
implies some compatibility conditions for $\zL$ and $\zG$ (cf.
(\ref{Jsc})). Poisson bracket are just Jacobi brackets with $\zG=0$,
i.e.  such that $\1$ is a central element, $\{\1,\cdot\}=0$.
The concepts of a Poisson and a Jacobi bracket can be easily
extended to the graded case.

\begin{definition}
A \emph{graded Jacobi bracket} of degree  $k$  on  a  $G$-graded
(think e.g.  $\Z$-graded) associative commutative algebra
$\A=\oplus_{g\in G}\A^g$ with unity $\1$ is a graded bilinear map
\be \{\cdot,\cdot\}:\A\ti\A\ra\A \ee of degree $k$, i.e.  $w(\{
a,b\})=w(a)+ w(b)+k$, such that
\begin{enumerate}
\item $\{ a,b\}=-(-1)^{\la w(a)+k,w(b)+k\ran}\{ b,a\}$ (graded
anticommutativity), \item $\{ a,bc\}=\{  a,b\}  c+(-1)^{\la
w(a)+k,b\ran}b\{  a,c\}-\{ a,\1\} bc$  (generalized graded
Leibniz rule), \item $\{\{ a,b\},c\}=\{ a,\{  b,c\}\}-(-1)^{\la
w(a)+k,w(b)+k\ran}\{ b,\{  a,c\}\}$ (graded Jacobi identity).
\end{enumerate}
Such a bracket is called a \emph{graded Poisson bracket} if $\1$ is its
central element. Note that $\Z$-graded
algebras  furnished  with  a  graded  Poisson bracket of degree
$-1$ are sometimes called {\em Gerstenhaber  algebras\/}  (see
\cite{KS1}, \cite{KS1a}).
\end{definition}

\begin{definition} If $H$ is an element in a graded Poisson algebra $\A$, we can consider the bracket $\{
a,b\}^H=\{\{ a,H\},b\}$, called the \emph{derived bracket}
(associated with the `Hamiltonian' $H$) \cite{KS}. The odd Hamiltonians we call
\emph{homological} if $\{ H,H\}=0$ (the latter condition is
nontrivial for odd Hamiltonians).
\end{definition}

\medskip\noindent{\bf Problem.}
Show that, if the graded Poisson bracket is of degree $k$ and $H$
is a homological Hamiltonian of degree $h$, then the derived
bracket is a graded Lie bracket of degree $h+2k$.


\begin{example} {\bf (The Schouten bracket)} Several natural graded Lie brackets of tensor fields
are associated with a given smooth ($C^\infty$) manifold $M$. Historically
the first one was probably the celebrated Schouten-Nijenhuis  bracket  $[\
,\ \rS$ defined on multivector  fields  (see  \cite{Ni,Sc}  for  the
original and \cite{Mi2} for a modern version). It is the unique
graded Poisson extension of the usual bracket $[\cdot  ,\cdot  ]$
of vector  fields  to  the  Grassmann  algebra
$\X(M)=\bigoplus_{n\in\N}\X^n(M)$  of  multivector   fields.
Consequently,

\begin{itemize}
\item the degree of $X\in \X^n(M)$ with respect to the  bracket
is $(n-1),$ \item
 $[X,f\rS  = X(f)$  for  $X\in  \X^1(M)$,   $f\in \X^0(M)=C^\infty  (M)$;
\item For $X\in  \X^k(M)$,  $Y\in \X^l(M)$, we have
\begin{equation}\label{Sch}[X,Y\wedge Z\rS =[X,Y\rS\wedge Z+(-1)^{(k-1)l}Y\wedge[X,Z\rS\,.
\end{equation}
\end{itemize}
Explicitly,
\begin{eqnarray}\nn
\lefteqn
{[X_1\wedge\ldots\wedge X_r,Y_1\wedge\cdots\wedge Y_n\rS =} \\
& &
\sum_{k,l}(-1)^{k+l}[X_k,Y_l]\wedge\ldots\wedge\widehat{X_k}\wedge
\ldots\we X_r\we Y_1\we\ldots\we\widehat{Y_l}\we\ldots\we Y_n,
\label{Schouten3}
\end{eqnarray}
where $X_k,Y_l\in\X^1(M)$ and `$\widehat{\ \ }$' stand for the omission. Note
that $[X,Y\rS$, with $X$ being a vector field, is just the Lie
derivative $\Ll_XY$.
\end{example}

It is easy to see that condition (\ref{ptnsr}) defining a
Poisson tensor can be rewritten in terms of the
Schouten-\-Nijenhuis bracket as
\begin{equation}\label{SP}[\zL,\zL]^{SN}=0\,.
\end{equation}
Moreover, the corresponding Poisson bracket can be viewed as the
derived bracket:
\begin{equation}\label{poder}\{ f,g\}_\zL=[[f,-\zL]^{SN},g]^{SN}\,.
\end{equation}
The Poisson tensor $-\zL$ plays the role of the homological
Hamiltonian which is quadratic (of degree 2). The derived bracket
is therefore of degree $2+2(-1)=0$, so closed on basic (degree 0)
functions.

Similarly to (\ref{SP}), (\ref{Jabra}) is a Jacobi bracket if and
only if
\begin{equation}\label{Jsc}
[\zL,\zL\rS=2\zL\we\zG\quad\text{and}\quad
[\zG,\zL\rS=\Ll_\zG\zL=0\,.
\end{equation}

\begin{remark} One can consider as well the \emph{symmetric Schouten bracket} (see \cite  {DV-M}).
It   is an   ordinary (non-graded) Lie bracket extending  the Lie
bracket of  vector  fields, defined on symmetric contravariant
tensors and   satisfying   an   analog   of (\ref{Schouten3}):
\begin{eqnarray}\nn
\lefteqn
{[X_1\vee\ldots\vee X_r,Y_1\vee\cdots\vee Y_n]^{SS} =} \\
& & \sum_{k,l}[X_k,Y_l]\vee\ldots\vee\widehat{X_k}\vee \ldots\we
X_r\we Y_1\we\ldots\we\widehat{Y_l}\we\ldots\we Y_n\,.
\label{SSchouten}
\end{eqnarray}
The symmetric Schouten bracket is, however, nothing but the
standard symplectic Poisson bracket on $\sT^\* M$ reduced to
polynomial functions. Polynomial functions on $\sT^\* M$
represent, namely, symmetric contravariant tensors by an extension
of (\ref{zi}),
\begin{equation}\label{zi1}\zi\left({f(q)\pa_{q^{a_1}}\vee\cdots\vee\pa_{q^{a_n}}}\right)=f(q)p_{a_1}\cdots p_{a_n}\,,
\end{equation}
and the brackets are identified according to (\ref{i1}) and
(\ref{i2}).
\end{remark}

\section{Algebroids}
\emph{Lie algebroids} are geometric objects which are so common and
natural that we are often working with them not even mentioning it. The
people told that they are using a Lie algebroid resemble Mr.
Jourdain who was surprised and delighted to learn that he has been
speaking prose all his life without knowing it. One can consider
also a more general object, a \emph{skew algebroid}, for which we drop
the Jacobi identity.

A \emph{Lie pseudoalgebra},  a pure  algebraic   counterpart of a Lie algebroid,
appeared first in the paper of Herz  \cite{He}
but one can find similar concepts under more than a dozen  of
names  in the   literature   (e.g.  Lie modules,  $(R,A)$-Lie
algebras, Lie-Cartan pairs, Lie-Rinehart algebras, differential
algebras, etc.). Lie algebroids were introduced  by Pradines
\cite{Pr1} as infinitesimal parts of differentiable groupoids. In
the same year the booklet \cite{Ne} by Nelson was published, where a
general theory of Lie modules together with a big part of the
corresponding differential calculus can be found. We also refer to
a survey article by Mackenzie \cite{Ma} and his book \cite{Mac2}.

\subsection{Skew algebroids}

Let $\zt:E\to M$ be a rank-$n$ vector bundle over an
$m$-dimensional manifold $M$ and let $\zp:E^\*\ra M$ be its dual.
Recall that the Grassmann algebra $\Gr(E)=\oplus_{i=0}^\infty\Sec(\we^i
E)$
of multisections of $E$ is a
graded commutative associative algebra with respect to the wedge
product.

There are different equivalent ways to define a \emph{skew
algebroid} structure on $E$. Here we will list only four of them.
The notation is borrowed from \cite{GG, GGU,GU2} and we refer to
these papers for details. In particular, we use affine coordinates
$(x^a,\zx_i)$ on $E^\*$ and the dual coordinates $(x^a,y^i)$ on
$E$, associated with dual local bases, $(e_i)$ and $(e^i)$, of
sections of $E$ and $E^\*$, respectively.

\begin{definition}  A \emph{skew algebroid} structure on $E$ is given by a linear bivector
field $\zP$ on $E^\*$. In local coordinates,
\begin{equation}\label{SA} \Pi =\frac{1}{2}c^k_{ij}(x)\zx_k
\partial _{\zx_i}\we \partial _{\zx_j} + \zr^b_i(x) \partial _{\zx_i}
\wedge \partial _{x^b}\,,
\end{equation} where $c^k_{ij}(x)=-c^k_{ji}(x)$. If \ $\Pi$ is a Poisson tensor, we speak about a \emph{Lie algebroid}.
\end{definition}

As the bivector field $\Pi$ defines a bilinear bracket
$\{\cdot,\cdot\}^\Pi$ on the algebra $C^\infty(E^\*)$ of smooth
functions on $E^\*$ by
$\{\zvf,\psi\}^{\zP}=\langle\zP,\xd\zvf\we\xd\psi\rangle$, where
$\langle\cdot,\cdot\rangle$ stands for the contraction, we get the
following.

\begin{theorem}\label{ty}
A skew algebroid structure $(E,\Pi)$ can be equivalently defined
as
\begin{itemize}
\item a skew-symmetric $\R$-bilinear bracket $[\cdot ,\cdot]^\Pi $
on the space $\Sec(E)$ of sections of $E$, together with a vector
bundle morphisms\ $\zr=\zr^\Pi \colon E\rightarrow T M$ (\emph{the
anchor}), such that
\begin{equation}\label{qd1} [X,fY]^\Pi =\zr^\Pi(X)(f)Y +f [X,Y]^\Pi\,,
\end{equation} for all $f \in C^\infty (M)$, $X,Y\in \Sec(E)$;

\item a graded skew-symmetric bracket $\lna\cdot ,\cdot\rna^\Pi$
of degree $-1$, the \emph{algebroid Schouten bracket}, on the
Grassmann algebra $\Gr(E)$, satisfying the Leibniz rule
\begin{equation}\label{SchL}\lna X,Y\wedge Z\rna^\zP =\lna X,Y\rna^\zP\wedge Z+(-1)^{(k-1)l}Y\wedge\lna X,Z\rna^\zP\,,
\end{equation}
for $X\in  \Gr^k(E)$,  $Y\in \Gr^l(E)$;

\item or as a graded derivation $\xdp$ of degree 1 in the
Grassmann algebra $\Gr(E^\*)$ (the \emph{de Rham derivative}),
\begin{equation}\label{dR}\xdp(\za\we\zb)=\xdp\za\we\zb+(-1)^{w(\za)}\za\we\xdp\zb\,.\end{equation}
\end{itemize}
Moreover, the following properties of the above structures are equivalent:
\begin{itemize}
\item $(E,\Pi)$ is a Lie algebroid. \item $[\cdot,\cdot]^\zP$ is a
Lie bracket. \item $\lna\cdot,\rna^\zP$ is a graded Poisson
bracket. \item $(\xdp)^2=0$.
\end{itemize}
In the latter case, the \emph{Lie algebroid cohomology} is defined in the
standar way: $H^\bullet(E;\xdp)=({\Ker\xdp}/{\im\xdp})^\bullet$.
\end{theorem}

The bracket $[\cdot,\cdot]^\Pi$ and the anchor $\zr^\Pi$ are related
to the bracket $\{\cdot,\cdot\}^{\zP}$ according to the formulae:
\begin{eqnarray}\label{drel}
        \zi([X,Y]^\Pi)&= \{\zi(X), \zi(Y)\}^{\zP},  \\
        \zp^\*(\zr^\Pi(X)(f))       &= \{\zi(X), \zp^\*f\}^{\zP}\,,\label{drel1}
                                                   \end{eqnarray}
where we denoted with $\zi(X)$ the linear function on $E^\*$
associated with the section $X$ of $E$, i.e.  $\zi(X)(e^\ast_p)=\la
X(p),e^\ast_p\ran$ for each $e^\ast_p\in E^\ast_p$.

The algebroid Schouten bracket is the unique graded extension of
$[\cdot,\cdot]^\zP$ satisfying the Leibniz rule.
The \emph{de Rham derivative} $\xdp$ is determined by the formula
\begin{equation}\label{dRd}(\xdp\zm)^v=[\zP,\zm^v]^{SN},
\end{equation}
where $\zm^v$ is the natural vertical lift of a `$k$-form'
$\zm\in\Gr^k(E^\*)$ to a vertical $k$-vector field on $E^\*$ and
$[\cdot,\cdot]^{SN}$ is the Schouten-Nijenhuis bracket of
multivector fields.
It can also be written in the \emph{Cartan form} \bea\nn
       & \xd^\zP\zm(X_1,\dots,X_{k+1}) = \sum_i (-1)^{i+1}
\zr^\Pi(X_i)(\zm(X_1,\dots,\widehat{X}_i,\dots ,X_{k+1})) \\
& + \sum_{i<j} (-1)^{i+j}\zm ([X_i,X_j]^\Pi, X_1,\dots ,
\widehat{X}_i, \dots ,\widehat{X}_j, \dots ,
X_{k+1})\,.\label{Cartan} \eea

\medskip
In local bases of sections and the corresponding local
coordinates,
\begin{eqnarray}\label{r1a}
[e_i,e_j]^\Pi(x)&=&c^k_{ij}(x)e_k,\\
\zr^\Pi(e_i)(x)&=&\zr^a_i(x)\pa_{x^a},\label{r2}\\
(\xdp f)(x)&=&\zr^a_i(x)\frac{\pa f}{\pa x^a}(x)e^i,\label{r3}\\
(\xdp e^i)(x)&=&c^i_{lk}(x)e^k\we e^l\,.\label{r4}
\end{eqnarray}

\medskip
Given a skew algebroid $E$, we can associate with any
$C^1$-function $H$ on $E^\*$ its \emph{Hamiltonian vector field}
$\X_H$ like in the standard case: $\X_H=i_{\xd\! H}\Pi$ that allows for a sort of `Hamiltonian mechanics'. In local coordinates,
\begin{equation}\label{ham}
\X_H(x,\zx)=\left(c^k_{ij}(x)\zx_k\frac{\pa H} {\partial
{\zx_i}}(x,\zx)- \zr^a_j(x)\frac{\pa H}{\partial
{x^a}}(x,\zx)\right) \partial _{\zx_j} + \zr^b_i(x) \frac{\pa
H}{\partial {\zx_i}}(x,\zx)
\partial _{x^b}\,.
\end{equation} Another geometrical construction in the skew-algebroid setting is the \emph{complete lift of an algebroid
section} (cf. \cite{GU1,GU2}). For every $C^1$-section,
$X=f^i(x)e_i\in\Sec(E)$, we can construct canonically a vector
field $\dd_T^\Pi(X)\in\Sec(T E)$ which in local coordinates reads
as
\begin{equation}\label{eqn:tan_lift}
\dd_\sT^\Pi(X)(x,y) = f^i(x)\rho^a_i(x)\pa_{x^a} + \left( y^i
\rho^a_i(x) \frac{\pa f^k}{\pa x^a}(x) + c^k_{ij}(x)y^if^j(x)
\right) \pa _{y^k}.
\end{equation}
The vector field $\dd_T^\Pi(X)$ is homogeneous (linear with
respect to $y$'s).

\begin{theorem} [\cite{GU1,GU2}]
The pair $(E,\zP)$ defines a Lie algebroid if and only if $\dd_\sT^\Pi([X,Y]^\Pi)
= [\dd_\sT^\Pi(X), \dd_\sT^\Pi(Y)]$ for all $X,Y\in \Sec(E)$.
\end{theorem}

\begin{example} Any tangent bundle $E =\sT M$ of a manifold $M$, with $\zr =\Id_{\sT M}$ and
the usual Lie bracket of vector fields, is a Lie algebroid.
\end{example}

\begin{example} Any Lie algebra, $E =\g$, considered as a vector bundle
over one point $M=\{ pt\}$ with the trivial anchor $\zr = 0$, is a
Lie algebroid. This Lie algebroid can be viewed as a reduction of
the tangent bundle of any Lie group $G$ associated with $\g$,
namely $\g=\sT G/G$, in which sections of $\g$ are interpreted as
invariant vector fields on $G$.
\end{example}

\begin{example} The above reduction procedure can be generalized to the case
of any principal $G$-bundle $P$. Invariant vector fields on $P$
are closed with respect to the Lie bracket and can be viewed as
sections of the vector bundle $E=\sT P/G$ which becomes a Lie
algebroid. This is the so called \emph{Atiyah algebroid}
associated with the principal bundle $P$.
\end{example}

\begin{example} There is a canonical Lie algebroid structure on the
cotangent bundle $\sT^*M$ associated with a Poisson tensor $\zL$ on
$M$. This is the unique Lie algebroid bracket $[\cdot,\cdot]_\zL$
of differential 1-forms for which $[\xd f,\xd g]_\zL=\xd\{
f,g\}_\zL$, where $\{\cdot,\cdot\}_\zL$ is the Poisson bracket of
functions for $\zL$ and the anchor map is $\widetilde{\zL}:\sT^*M\ra
\sT M$. Explicitly (cf. (\ref{Kb})),
\be
[\za,\zb]_\zL=\Ll_{\widetilde{\zL}(\za)}\zb-\Ll_{\widetilde{\zL}(\zb)}\za-\xd\la\zL,\za\we\zb\ran\,.
\end{equation} This Lie bracket was defined first by Fuchssteiner \cite{Fu} but it is usually called the \emph{Koszul bracket} \cite{Ko}. The corresponding linear Poisson tensor on $\sT M$ is the tangent lift $\dd_\sT\zL$ \cite{GU1,GU2}. As the tangent lift respects the Schouten bracket \cite{G-U1}, it is again a Poisson structure, this time linear.
\end{example}

In the next sections we will show two more interpretations of a skew
algebroid: as a morphism of \emph{double vector bundles} and as a
vector field on a \emph{graded manifold}.

\subsection{Differential calculus on Lie algebroids}

        Let  us consider a Lie algebroid structure $(E,\zP)$ on a vector bundle $\zt\colon E\rightarrow M$, with the Lie bracket $[\cdot,\cdot ]^\zP$ on sections of $E$ and the anchor
$\zr^\zP \colon E\rightarrow \sT M$.  Then, we can construct a
well-known generalization of the standard Cartan calculus of
differential forms and vector fields (see e.g. \cite{Mac2,Mr1}).

First, we have the exterior (de Rham) derivative $\xd^\zP \colon
\Gr^k(E^\*) \rightarrow \Gr^{k+1}(E^\*)$ (\ref{Cartan}).  For $X
\in \Gr^k(E)$, the contraction $\xi_X \colon  \Gr^p(E^\*)
\rightarrow \Gr^{p-k}(E^\*)$ is defined in the standard way for $k=1$ and extended by $\xi_{X_1\we\cdots\we X_k}=\xi_{X_1}\cdots \xi_{X_k}$ for $X_i\in\Sec(E)$ (this produces a sign factor with respect to another convention for the contraction).
The Lie differential operator
                \be \Ll^\zP_X\colon \Gr^p(E^\*)\rightarrow \Gr^{p-k+1}(E^\*)
                                \ee
         is defined as the graded commutator
                \be \Ll^\zP_X = [\xi_X,\xd^\zP]= \xi_X \circ \xd^\zP - (-1)^k \xd^\zP \circ
\xi_X.
                                \ee
        The following proposition contains a list of well-known properties of these objects.

\begin{proposition}
        Let  $\zm\in \Gr^k(E^\*),\ \zn\in \Gr(E^\*)$ and $X,Y \in
\Gr^1(E)$. We have
        \begin{enumerate}
        \item $\xd^\zP \circ \xd^\zP =0$\,,
        \item $\xd^\zP(\zm\wedge \zn) = \xd^\zP\zm \wedge \zn + (-1)^k\zm
\wedge \xd^\zP\zn$\,,
        \item $\xi_X(\zm\wedge \zn) = \xi_X\zm \wedge \zn + (-1)^k\zm
\wedge \xi_X\zn$\,,
        \item $\Ll^\zP_X(\zm\wedge \zn) = \Ll^\zP_X\zm\wedge \zn +
\zm\wedge \Ll^\zP_X \zn$\,,
        \item $[\Ll^\zP_X, \Ll^\zP_Y]=\Ll^\zP_X\circ \Ll^\zP_Y - \Ll^\zP_Y \circ \Ll^\zP_X =\Ll^\zP_{[X,Y]^\zP}$\,,
        \item $[\Ll^\zP_X, \xi_Y]=\Ll^\zP_X \circ \xi_Y - \xi_Y\circ \Ll^\zP_X = \xi_{[X,Y]^\zP}$\,.
        \end{enumerate}
                 \end{proposition}

       \noindent The last formula can be generalized in the following way
(cf. \cite{G-U2,Mr1,Mi2}).

        \begin{theorem}\label{tsb}
        For $X\in \Gr^{k+1}(E) $ and $Y\in \Gr^{l+1}(E)$,
                \be  [\Ll^\zP_X, \xi_Y]=\Ll^\zP_X\circ \xi_Y   - (-1)^{(k+1)l}\xi_Y\circ \Ll^\zP_X =
 \xi_{\lna X,Y\rna^\zP}\,,
                                \ee
                                where $\lna\cdot,\cdot\rna^\zP$ is the algebroid Schouten bracket. In particular, for $X\in \Gr^1(E)$ and $f\in
\Gr^0(E) = C^\infty(M)$  we have
                \be \lna X,f\rna^\zP = \zr^\Pi(X)(f)\,.
                                \ee
 \end{theorem}

        There is also a \emph{symmetric Schouten bracket} which extends
the Schouten bracket on $ \Gr^0(E)\oplus \Gr^1(E)$ to symmetric
multisections. This bracket is just the polynomial part of the
Poisson bracket $\{\cdot,\cdot\}^\zP$.
Before we pass to other brackets, let us introduce the bi-graded space \be\zF(E)=\oplus_{k,l=0}^\infty\zF^k_l(E)=\Gr(E\oplus_M E^\*)\,, \quad
\zF^k_l(E) = \Sec(\wedge^kE\otimes_M \wedge^lE^\*)\,,
\ee
of tensor fields of mixed type. Of course, we can identify $\zF^k_l(E)$ and
$\zF^l_k(E^\*)$. For   $K\in \zF^k_1(E^\*)$, we define
the contraction
                \be \xi_K\colon \Gr^n(E^\*) \rightarrow \Gr^{n+k-1}(E^\*)
                                        \ee
       in a natural way: for simple tensors  $K = \zm\otimes X$, where $\zm \in \Gr^k
(E^\*),\, X\in \Gr^1(E) $, we just put
                \be \xi_K\zn = \zm\wedge \xi_X\zn .
                         \label{1.8}\end{equation}         The corresponding Lie differential is defined by the formula
                \be \Ll^\zP_K = \xi_K \circ \xd^\zP + (-1)^k\xd^\zP \circ \xi_K
                                        \ee
        and, in particular,
                \be \Ll^\zP_{\zm\otimes X  } = \zm\wedge \Ll^\zP_X + (-1)^k \xd^\zP
\zm \wedge \xi_X.
                        \label{1.9}\end{equation}         This definition is compatible with the previous one in the
case of $K\in \zF^0_1(E^\*) = \Gr^1(E)=\Sec(E)$.
The contraction (insertion) $\xi_K$ can be extended to an operator
                \be \xi_K\colon \zF^n_1(E^\*) \rightarrow \zF_1^{n+k-1}(E^\*)
                                        \ee
        by the formula
                \be \xi_K(\zm\otimes X) = \xi_K(\zm)\otimes X .
                                        \ee

        \begin{theorem}
        The bracket
                \be [\cdot,\cdot]^{NR} \colon \zF_1^{k+1}(E^\*) \times
\zF_1^{l+1}(E^\*) \rightarrow \zF_1^{k+l+1}(E^\*),
                                        \ee
given by the formula
                \be [K,L]^{NR} = \xi_KL- (-1)^{kl}\xi_LK,
                         \label{1.10}\end{equation} defines a graded Lie algebra structure on the graded space $\zF_1(E^\*) = \oplus_{k\in \N}
\zF_1^{k}(E^\*)$. For simple tensors $\zm\otimes X \in
\zF_1^{k}(E^\*)$ and $\zn \otimes Y \in \zF_1^{l}(E^\*)$, we get
                \be[\zm\otimes X, \zn\otimes Y]^{NR} = \zm\wedge
\xi_X\zn \otimes  Y + (-1)^k \xi_Y \zm\wedge \zn \otimes X.
                         \label{1.11}\end{equation}         \end{theorem}
The bracket $[\cdot,\cdot]^{NR}$ is called  the \emph{(generalized)
Nijenhuis-Richardson bracket}.

\begin{remark} The generalized  Nijenhuis-Richardson bracket is a purely vector bundle bracket and does not
depend on the Lie algebroid structure. It is a geometric
counterpart of the purely algebraic bracket (\ref{aNR}). For $E
=\sT M$, we get the classical Nijenhuis-Richardson bracket of
vector-valued forms \cite{N-R}.
\end{remark}

        Another important bracket, the  \emph{generalized Fr\"olicher-Nijenhuis bracket}, is also a bracket on the graded space $\zF_1^{}(E^\*) =\oplus_{k\in \N}\zF_1^{k}(E^\*)$ of `vector-valued forms',
defined for simple tensors $\zm \otimes X \in \zF_1^{k}(E^\*)$ and
$\zn \otimes Y \in \zF_1^{l}(E^\*)$ by \bea\nn
       & [\zm \otimes X, \zn \otimes Y]^{FN}_\zP =
        \zm\wedge \zn \otimes
[X,Y]^\zP + \zm\wedge \Ll^\zP_X\zn\otimes Y -\Ll^\zP_Y\zm \wedge
\zn \otimes X\\& + (-1)^k(\xd^\zP\zm \wedge \xi_X\zn \otimes Y +
\xi_Y\zm
\wedge \xd^\zP\zn \otimes X) \nn \\
        &= (\Ll^\zP_{\zm\otimes X}\zn) \otimes Y - (-1)^{kl}(\Ll^\zP_{\zn
\otimes Y}\zm)\otimes X + \zm\wedge \zn \otimes [X,Y]^\zP.\label{88}
                                    \eea

        \begin{theorem} [\cite{GU2,K-M-S}]
        The  formula (\ref{88}) defines a graded Lie bracket
of degree 0 on the graded space $\zF_1^{}(E^\*)
=\oplus_{k\in \N}\zF_1^{k}(E^\*)$ of vector-valued forms. Moreover,
\begin{eqnarray} [\Ll^\zP_K,\Ll^\zP_L]& = &  \Ll^\zP_K\circ\Ll^\zP_L-(-1)^{kl}\Ll^\zP_L
\circ\Ll^\zP_K= \nn \\
&=& \Ll^\zP_{[K,L\rF_\zP}\,, \\ \label{Fro1} [\Ll^\zP_K,i_L] & = &
\Ll^\zP_K\circ i_L-(-1)^{k(l+1)}i_L
\circ\Ll^\zP_K\nn \\
&=& i_{[K,L\rF_\zP }-(-1)^{k(l+1)}\Ll^\zP_{i_LK}\,.\label{Fro2}
\end{eqnarray}
\end{theorem}

\noindent{\bf Problem.} Prove that, for $N$ being a $(1,1)$ tensor
interpreted as a morphism $N:E\to E$ of vector bundles, we have
\begin{equation}\label{Nij} [N,N]^{FN}_\zP(X,Y)=[NX,NY]^\zP-N([NX,Y]^\zP+[X,NY]^\zP-N[X,Y]^\zP)\,,
\end{equation}
for any $X,Y\in\Sec(E)$. The tensor $[N,N]^{FN}_\zP$ is sometimes called the \emph{(generalized) Nijenhius torsion} of $N$.

\medskip
In the case of the canonical Lie algebroid $E=\sT M$, we obtain the classical Fr\"olicher -Nijenhuis bracket on
the graded space $\zW(M;\sT M)=\zF_1(\sT^\* M)$ of vector-valued
forms \cite{DV-M,F-N,K-M-S,Mi1}.

Note that there are some interesting relations of the classical
Nijenhuis-Richardson and Fr\"olicher-Nijenhuis brackets on $M$
with the Schouten bracket on $\sT^\* M$.
Let us first recall that
any vector field $X$ on $M$ can be identified with a linear function $\zi(X)$ on $\sT^\* M$. As $\sT^\* M$ is
canonically a symplectic manifold, we can associate with $\zi(X)$
its Hamiltonian vector field which will be denoted $\H(X)$ and
called the \emph{cotangent lift} of $X$.
Second, any one-form $\za$ can be lifted to a vertical vector
field $\V(\za)$, \be\V(f_a(q)\xd q^a)=f_a(q)\pa_{p_a}\,.\ee This
vertical lift can be extended to $k$-forms by
\be\V(\za_1\we\cdots\we\za_k)=\V(\za_1)\we\cdots\we \V(\za_k)\,.\ee
We can extend the maps $\zi$ and $\H$ to linear maps
$\J,\H:\zW(M;TM)\to\X(\sT^\* M)$ by
                \be \J(\zm \otimes X) = - \zi(X)\V(\zm)
                                        \label{6.3}\ee
        and
                \be \H(\zm \otimes X)  = \H(X)
\wedge \V(\zm) - \zi(X) \V (\xd \zm),
                                        \label{6.4}\ee
        for simple tensors $\zm \otimes X \in \zW(M;TM)$.

\begin{theorem}[\cite{G-U2}]
        The mappings $\J,\H:\zW(M;TM)\to\X(\sT^\* M)$
        are injective homomorphisms (embeddings) of, respectively, the
        Nijenhuis-Richardson and the Fr\"olicher-Nijenhuis
bracket into the Schouten-Nijenhuis bra\-cket:
                \begin{eqnarray}  \J([K, L]^{NR}) &=& [\J(K), \J(L)]^{SN}\,,\\
                  \H([K, L]^{FN}) &=& [\H(K), \H(L)]^{SN}\,.
                                        \end{eqnarray} \end{theorem}

\subsection{Nijenhuis tensors}

If, for a Lie algebroid $(E,\zP)$, the Nijenhuis torsion
(\ref{Nij}) of a tensor $N:E\to E$ vanishes, we call $N$ a
\emph{Nijenhuis tensor} (see \cite{KSM}). The crucial property of a
Nijenhuis tensor as defining a contraction of the bracket is the following (cf. \cite{KSM,GU1}).

\begin{theorem} If $N$ is a Nijenhuis tensor for a Lie algebroid bracket
$[\cdot,\cdot]^\zP$ on $E$ with an anchor map $\zr^\Pi:E\ra \sT M$,
then the contracted bracket
\begin{equation}\label{Nij cont}[X,Y]^\zP_N=[NX,Y]^\zP+[X,NY]^\zP-N[X,Y]^\zP
\end{equation}
is again a Lie algebroid bracket on $E$ with the anchor
$\zr^\Pi_N=\zr^\Pi\circ N$. This bracket corresponds to the linear Poisson structure
$\zP_N=\Ll_{\J(N)}\zP$. Moreover, $N:E\to E$ is a morphism of the
Lie algebroid $(E,\zP_N)$ into the Lie algebroid $(E,\zP)$:
\be
[NX,NY]^\zP=N\left([X,Y]^\zP_N\right)\,.
\ee
\end{theorem}

\begin{remark} If the contracted bracket $[\cdot,\cdot]^\zP_N$ is again
a Lie bracket, then $N$ is called \emph{weak Niejnhuis} (cf. \cite{CGMc}).
The above theorem implies that Nijenhuis tensors are weak-Nijenhuis.
Tensors $N:\sT M\to\sT M$ satisfying $N^2=-Id$ are called
\emph{almost complex structures}. The celebrated
Newlander-Nirenberg theorem states that an almost complex
structure $N$ is integrable, i.e.  comes from a true complex structure,
if and only if $N$ is Nijenhuis.
\end{remark}

\begin{example} {\bf (Frobenius manifolds)} An algebraical part of the structure of a {\it Frobenius
manifold} consists of a unital commutative associative
multiplication "$\circ$" in the space $\X^1(M)$ of vector fields
which comes from a symmetric vector valued two-form
$C\in\Sec(\vee^2\sT^\*M\ot_M\sT M)$. This multiplication is supposed
to satisfy the following axiom proposed by  Hertling and Manin:
\begin{equation}\label{Faxiom} \Ll_{X\circ Y}C = X\circ\Ll_Y C+\Ll_XC\circ
Y\,.
\end{equation}  In terms of structure functions in local coordinates, (\ref{Faxiom}) reads as
\be\label{YA}
C^m_{sj}\frac{\pa C^s_{lr}}{\pa x^k}+
C^m_{sk}\frac{\pa C^s_{lr}}{\pa x^j} - C^m_{sr}\frac{\pa
C^s_{jk}}{\pa x^l} -C^m_{sl} \frac{\pa C^s_{jk}}{\pa x^r} +
\frac{\pa C^m_{jk}}{\pa x^s}C^s_{lr} - \frac{C^m_{lr}}{\pa
x^s}C^s_{jk}=0\,,
\ee
for all $m,j,k,l,r$. An interpretation of
the above conditions can be found in an old paper by Yano and Ako
\cite{YA}, where they constructed several classes of "differential
concomitants" in the sense of Schouten. One among them leads
exactly to (\ref{YA}), so that the above differential constraints
on the tensor field $C$ are sometimes referred to as the {\it
Yano-Ako conditions} (see \cite{Mag1} and the discussion there). We
will show that (\ref{Faxiom}) can be interpreted as
vanishing of a Nijenhuis torsion.

Recall first that symmetric multi-vector fields on the manifold $M$
with the symmetric Schouten bracket  can be identified with the
graded algebra $\A=\oplus_{k=0}^\infty\A^k$ of polynomial
functions on $\sT^\* M$ with the symplectic Poisson bracket.
Second, the unital commutative associative
multiplication $C$ in $\X^1(M)=\A^1$ defines an
$\A^0=C^\infty(M)$-linear projection $N=N(C):\A\to\A^1$ defined
by:
\begin{equation}\label{proj} N(1)=E_0\,,\qquad N(X_1\cdots X_k)=X_1\circ\cdots\circ
X_k\quad\text{for all}\quad X_1,\dots,X_k\in\A^1\,,
\end{equation} where $E_0\in\A^1$ corresponds to the unity vector field of the multiplication. As $\A^1$ is a Lie subalgebra
in $\A$ with the Poisson bracket, $N$ is a Nijenhuis tensor for
the Lie algebra structure in $\A$ if and only if $\Ker(N(C))$ is
also a Lie subalgebra \cite{CGM}. It can be directly checked that the latter
is equivalent to (\ref{Faxiom}).

\begin{theorem}
$N(C)$ is a Nijenhuis tensor for the Lie algebra
$(\A,\{\cdot,\cdot\})$ if and only if $C$ satisfies  the Yano-Ako conditions
(\ref{Faxiom}).
\end{theorem}

Note that the above observation is closely related to the so
called {\it coisotropic deformations} of associative structures as
studied e.g. in \cite{KoMa} and that one can easily prove also quantum or
supersymmetric analogs of the above theorem.  More about
Nijenhuis tensors for general brackets can be found in
\cite{CGM, KSM}. Note finally that Nijenhuis tensors compatible
with Poisson structures, \emph{Poisson-Nijenhuis tensors}, provide
a useful language for studying integrability of Hamiltonian
systems \cite{KSM,MM}.
\end{example}

\section{Double vector bundles and formalisms of Mechanics}

The starting point of what follows is the observation that a skew (or Lie) algebroid can be described as a particular morphism of double vector bundles.

\begin{definition} A \emph{double vector bundle} is a manifold $K$
with two compatible vector bundle structures, $\zt_i:K\to K_i$, $i=1,2$.
The compatibility means that the Euler (Liouville) vector fields
(generators of homotheties) associated with the two vector bundle structures
commute.
\end{definition}

This definition implies that, with every double vector bundle, we
can associate the following diagram of vector bundles in which
both pairs of parallel arrows form vector bundle morphisms:

\begin{equation}\label{F1.3x}\xymatrix{
 & K\ar[dr]^{\tau_2}\ar[dl]_{\tau_1} & \\
K_1\ar[dr]^{\tau'_2} &  & K_2\ar[dl]_{\tau'_1} \\
 & M & }
\end{equation}

The above geometric definition (cf. \cite{GR,GR1}) is a simplification of
the original categorical concept of a double vector bundle due to
Pradines \cite{Pr1}, see also \cite{Ma,KU}.

\begin{example}
Let $M$ be a smooth manifold and let $(x^a), \ a=1,\dots,m$, be a
coordinate system in $M$. We denote by $\zt_M \colon \sT M
\rightarrow M$ the tangent vector bundle and by $\zp_M \colon
\sT^\* M\rightarrow M$ the cotangent vector bundle. We have the
induced (adapted) coordinate systems $(x^a, {\dot x}^b)$ in $\sT
M$ and $(x^a, p_b)$ in $\sT^\* M$.
        Let $\zt\colon E \rightarrow M$ be a vector bundle and let $\zp
\colon E^\* \rightarrow M$ be the dual bundle.
  Let $(e_1,\dots,e_n)$  be a basis of local sections of $\zt\colon E\rightarrow M$ and let $(e^{1},\dots, e^{n})$ be the dual
basis of local sections of $\zp\colon E^\*\rightarrow M$. We have
the induced coordinate systems:
    \beas
    (x^a, y^i),\quad & y^i=\zi(e^{i})\,, \quad \text{in} \ E\,,\\
    (x^a, \zx_i), \quad &\zx_i = \zi(e_i)\,,\quad \text{in} \ E^\* \,.
    \eeas
    Thus we have the adapted local coordinates
    \beas
    (x^a, y^i,{\dot x}^b, {\dot y}^j ) &  \quad \text{in} \ \sT E\,,\\
    (x^a, \zx_i, {\dot x}^b, {\dot \zx}_j) & \quad \text{in} \ \sT E^\* \,,\\
    (x^a, y^i, p_b, \zp_j) & \quad \text{in}\ \sT^\*E\,,\\
    (x^a, \zx_i, p_b, \zf^j) & \quad \text{in}\ \sT^\* E^\* \,.
    \eeas
It is well known (cf. \cite{KT,KU,Ur}) that the tangent   bundle
$\sT E$ and the cotangent bundle $\sT^\*E$ are canonical examples
of double vector bundles: \be\xymatrix{
\sT E\ar[rr]^{\zt_E} \ar[d]_{\sT\zt} && E\ar[d]^{\zt} \\
\sT M\ar[rr]^{\zt_M} && M } \qquad {,}\qquad \xymatrix{
\sT^\ast E\ar[rr]^{\sT^\ast\zt} \ar[d]_{\zp_{E}} && E^\ast\ar[d]^{\zp} \\
E\ar[rr]^{\zt} && M } \ee with projections \be\zt_E(x^a, y^i,{\dot
x}^b, {\dot y}^j )=(x^a,y^i)\,,\quad \sT\zt(x^a, y^i,{\dot x}^b,
{\dot y}^j )=(x^a,{\dot x}^b)\ee and \be\zp_{E}(x^a, y^i, p_b,
\zp_j)=(x^a, y^i)\,,\quad \sT^\ast\zt(x^a, y^i, p_b,
\zp_j)=(x^a,\zp_j)\,.\ee The corresponding pairs of commuting
Euler vector fields are, respectively,
\be\nabla_1=\dot
x^a\partial_{\dot x^a}+\dot y^i\partial_{\dot y^i}\,,\quad
\nabla_2= y^i\partial_{y^i}+\dot y^j\partial_{\dot y^j}\,,\ee
and
\be\nabla'_1= p_a\partial_{p_a}+ \zp_i\partial_{ \zp_i}\,,\quad
\nabla'_2= p_a\partial_{p_a}+y^i\partial_{y^i}\,.\ee
\end{example}

The fundamental fact we will explore is that the double vector bundles $\sT^\ast
E^\ast$ and $\sT^\ast E$ are canonically isomorphic with an
isomorphism
 \be\label{iso}\cR_\zt \colon \sT^\*E \longrightarrow \sT^\* E^\*
                    \ee
  being simultaneously an anti-symplectomorphism  (we can choose a symplectomorphism as well) \cite{KT,KU,GU2}.
  In local coordinates, $\cR_\zt$ is given by
    \be\cR_\zt(x^a, y^i, p_b, \zp_j) = (x^a, \zp_i, -p_b,y^j).
                              \ee
This means that we can identify coordinates $\zp_j$ with $\zx_j$,
coordinates $\zf^j$ with $y^j$, and use the coordinates $(x^a,
y^i, p_b, \zx_j)$ in $\sT^\ast E$ and the coordinates $(x^a,
\zx_i, p_b,y^j)$ in $\sT^\ast E^\ast$, in the full agreement with
(\ref{iso}).

We known that skew algebroid structures on the vector bundle $E$
correspond to linear bivector fields on $E^\*$. As a matter of
fact, a 2-contravariant tensor $\zP$ on $E^\*$ is \emph{linear} if
and only if the corresponding mapping $\widetilde{\zP} \colon
\sT^\* E^\* \rightarrow \sT E^\*$ is a morphism of double vector
bundles. The commutative diagram \be\xymatrix{
\sT^\ast E^\ast\ar[r]^{\widetilde\Pi} \ar[d]_{\cR_\tau} & \sT E^\ast \\
\sT^\ast E\ar[ur]^{\ze} & } \ee describes a one-to-one
correspondence between linear 2-contravariant tensors $\zP_\ze$ on
$E^\*$ and homomorphisms  $\ze:\sT^\*E\to \sT E^\*$ of double vector bundles covering the
identity on $E^\*$ (cf. \cite{KU, GU2}).
In local coordinates, every  $\ze$ as above is of the form
                \be (x^a, \zx_i, {\dot x}^b, {\dot \zx}_j) \circ \ze = (x^a,
\zp_i, \zr^b_k(x)y^k, c^k_{ij}(x) y^i\zp_k + \zs^a_j(x) p_a)
                                                             \label{F1.4}\end{equation}
which shows that it covers also $\zr:E\to\sT M$ and  corresponds to the linear tensor
                \be \zP_\ze = c^k_{ij}(x)\zx_k \partial _{\zx_i}\otimes \partial
_{\zx_j} + \zr^b_i(x) \partial _{\zx_i} \otimes \partial _{x^b} -
\zs^a_j(x) \partial _{x^a} \otimes \partial _{\zx_j}
                                                             \label{F1.5}\end{equation}
 on $E^\*$. In \cite{GU2}, a (general) \emph{algebroid} is defined as the above morphism $\ze$ of
double vector bundles covering the identity on $E^\*$, while
a \emph{skew algebroid} (resp., \emph{Lie algebroid}) is such an
algebroid for which the tensor $\zP_\ze$ is skew-symmetric
(resp., Poisson).

\subsection{Lagrangian formalism for general algebroids}
A generalized Lagrangian formalisms for Lie algebroids has been proposed by Liberman and
Weinstein \cite{Lib,We} and developed in this setting by many
authors (e.g. \cite{LMM, Mar1, Mar2}).
In \cite{GG,GGU}, in turn, has been observed that its geometric background is actually based on double vector bundle morphisms $\ze$
and the Jacobi identity plays no role in the construction of dynamics, that gives a space for further generalizations.

\medskip
For a given an algebroid associated with the morphism $\ze:\sT^\*E\to \sT E^\*$, a Lagrangian $L:E\ra{\R}$ defines two smooth maps: the {\em
Legendre mapping}: { $\lambda_L:E\longrightarrow E^\ast$,
$\lambda_L=\tau_{E^\ast}\circ\ze\circ\xd L$,} and
the {\em Tulczyjew differential} { $\Lambda_{L}: E\longrightarrow
\sT E^\ast$, $\Lambda_{L}=\ze\circ\xd L$}. On the diagram it looks like
\begin{equation}\label{diag}\xymatrix{
\sT^\ast E\ar[rr]^{\ze}  && \sT E^\ast\ar[d]^{\zt_{E^\ast}} \\
E\ar[rr]^{\zl_L}\ar[u]^{\xd L}\ar@{.>}[rru]^{\Lambda_{L}} &&
E^\ast }.
\end{equation} The lagrangian function $L$ defines {\em the
phase dynamics}  as the set { $\cD=\Lambda_{L}(E)\subset \sT E^\ast$} which
can be understood as an implicit differential equation on $E^\*$,
solutions of which are `phase trajectories' of the system,
$\zb:\R\ra E^\ast$, and satisfy $\st(\zb)(t)\in \cD$, where $\st$ denotes the tangent prolongation of a $C^1$-curve. An analog of
the Euler-Lagrange equation for curves {\ $\gamma: \R\rightarrow
E$} is then
\be(E_L):\qquad
\st(\zl_L\circ\gamma)=\Lambda_L\circ\gamma\,.\ee
Equation
$(E_L)$ simply means that $\zL_L\circ\zg$ is an admissible curve
in $\sT E^\ast$, thus it is the tangent prolongation of
$\zl_L\circ\zg$. In local coordinates, $\cD$ has a
parametrization by $(x^a,y^k)$ via $\zL_L$ in the form (cf.
(\ref{F1.4}))
\begin{equation}\zL_L(x^a,y^i)= \left(x^a,\frac{\partial L}{\partial y^i}(x,y),
\zr^b_k(x)y^k, c^k_{ij}(x) y^i\frac{\partial L}{\partial y^k}(x,y)
+ \zs^a_j(x)\frac{\partial L}{\partial x^a}(x,y)\right)
\label{F1.4a}\end{equation} and equation $(E_L)$, for
$\zg(t)=(x^a(t),y^i(t))$, reads \be (E_L):\qquad\frac{\xd x^a}{\xd
t}=\zr^a_k(x)y^k, \quad \frac{\xd}{\xd t}\left(\frac{\partial
L}{\partial y^j}\right)= c^k_{ij}(x) y^i\frac{\partial L}{\partial
y^k} + \zs^a_j(x)\frac{\partial L}{\partial
x^a}\,.\label{EL2}\end{equation}
As one can easily see from (\ref{EL2}),
solutions are automatically admissible curves in $E$, i.e. 
$\zr^\Pi(\zg(t))=\st(\zt\circ\zg)(t)$. Since a curve in the canonical Lie
algebroid $\sT M$ is admissible if and only if it is a tangent
prolongation of its projection to $M$, first-order differential
equations for admissible curves in $\sT M$ may be viewed
as certain second-order differential equations for curves
in $M$. This explains why, classically, the Euler-Lagrange
equations are usually viewed as second-order equations.

\medskip\noindent
{\bf Remark.} In the standard case, $E=\sT M$, the Tulczyjew differential $\zL_L:\sT
M\ra\sT\sT^\ast M$ is sometimes called the
\emph{time evolution operator $K$} (see \cite{BGPR}), as the first
ideas of this operator go back to a work by Kamimura. This
operator has been studied by several authors in many variational
contexts, however, without recognition of its direct relation to a
(Lie) algebroid structure. We named this map after
Tulczyjew, since the above picture of the Lagrangian formalism is based on his ideas
\cite{Tu3}.

\begin{example}\label{e1} There are many examples based on Lie algebroids, e.g.
\cite{LMM,IMMS, Mar1}.
In particular, for the canonical Lie algebroid and the corresponding morphism which is the inverse of the Tulczyjew isomorphism \cite{Tu1}
\be\ze=\za_M^{-1}:\sT^\ast\sT
M\ra\sT\sT^\ast M\,,\ee
with $y^a=\dot{x}^a$, we get the
traditional Euler-Lagrange equations \be\frac{\xd x^a}{\xd
t}=\dot{x}^a, \quad \frac{\xd}{\xd t}\left(\frac{\partial
L}{\partial \dot{x}^a}\right)=\frac{\partial L}{\partial
x^a}\,.\ee

\medskip\noindent For a Lie algebroid which is just a Lie algebra with structure constants $c^k_{ij}$ with respect to a chosen
basis, we get the Euler-Poincar\'e equations \be\frac{\xd}{\xd
t}\left(\frac{\partial L}{\partial y^j}\right)= c^k_{ij}
y^i\frac{\partial L}{\partial y^k}\,. \ee
\end{example}

\medskip
The above examples are associated with Lie algebroids, but the presence of some
"nonholonomic constraints"  may lead to
Lagrangian systems on skew algebroids which are not Lie. This is related to
`quasi-Poisson brackets' associated with nonholonomic constraints
\cite{Mr,SM}.

\begin{example}{\bf (Skew algebroid of linear constraints)}  Consider an algebroid structure on a vector bundle $E$ equipped
with a Riemannian metric $\la\cdot,\cdot\ran_E$ and a vector
subbundle $C$ of $E$ (linear constraints). Let $P:E\ra C$ be the orthogonal projection.
We can choose a local basis of orthonormal sections
$(e_i)=(e_\za,e_A)$ of $E$ such that $(e_\za)$ is a basis of local
sections of $C$ and, identifying $E$ with $E^\*$, consider the corresponding affine coordinates $(x^a,y^k)=(x^a,y^\za,y^A)$ on $E$. According to the \emph{d'Alembert principle}, $\zd
L(\st(\zg)(t))\in C^0$, where $C^0\subset E^\ast$ is the
annihilator of $C$, the constrained dynamics is locally written
(cf. (\ref{EL2})) as
\be y^A=0\,,\quad\frac{\xd x^a}{\xd
t}=\zr^a_\za(x)y^\za, \quad \frac{\xd}{\xd t}\left(\frac{\partial
L}{\partial y^\zb}\right)- c^k_{\za\zb}(x) y^\za\frac{\partial
L}{\partial y^k} - \zs^a_\zb(x)\frac{\partial L}{\partial
x^a}=0\,.\label{EL3}\end{equation}
If we deal with a Lagrangian of "mechanical type",
\be
L=\frac{1}{2}(y^i)^2-V(x)\,,\ee
then $\frac{\partial L}{\partial
y^A}=y^A=0$ and equations (\ref{EL3}) reduce to
\be
y^A=0\,,\quad\frac{\xd x^a}{\xd t}=\zr^a_\za(x)y^\za\,, \quad \frac{\xd y^\zb}{\xd
t}=c^\zg_{\za\zb}(x) y^\za y^\zg -
\zs^a_\zb(x)\frac{\partial V}{\partial x^a}\,,\ee
that can be
viewed as Euler-Lagrange equations for the algebroid associated
with the orthogonal projection of the tensor $\zP_\ze$ onto
$C^\ast$, according to the orthogonal decomposition
$E^\ast=C^0\oplus C^\ast$ \cite{GG}. Of course, even for $E$ being a Lie
algebroid, if $C$ is not a Lie subalgebroid, the the projected tensor
is not a Poisson tensor and we deal with mechanics on a general
algebroid, in fact, a skew algebroid, since the
projected Poisson tensor remains skew-symmetric.
\end{example}

\subsection{Hamiltonian formalism for general algebroids}
Note that the linear tensor $\zP_\ze$ on $E^\*$ gives rise
also to a kind of Hamiltonian formalism. In \cite{GU2} and
\cite{OPB} one refers to a 2-contravariant tensor as to a
\emph{Leibniz structure}, that however may cause some confusion
with the \emph{Leibniz algebra} in the sense of Loday. Anyhow, in the
presence of $\zP_\ze$, by the \emph{Hamiltonian vector field}
associated with a function $H$ on $E^\*$ we understand the
contraction $\ix_{\xd H}\zP_\ze$ as in (\ref{ham}). Thus the
question of the Hamiltonian description of the dynamics $\cD\subset\sT E^\ast$ is the
question if $\cD$ is the image of a Hamiltonian vector field, i.e. 
\begin{equation}\mathcal{D}=\widetilde{\zP}_\ze(\xd H(E^\*))\,.\end{equation}
Every such a function $H$ we call a \emph{Hamiltonian  associated
with the Lagrangian} $L$. However, it should be stressed that,
since $\ze$ and $\zP_\ze$ can be degenerate, we have much more
freedom in choosing generating objects (Lagrangians and
Hamiltonians) than in the symplectic case. For instance, the
Hamiltonian is defined not up to a constant but up to a Casimir
function of the tensor $\zP_\ze$ and for the choice of the
Lagrangian we have a similar freedom. However, in the case of a
hyperregular Lagrangian, we recover the standard correspondence
between Lagrangians and Hamiltonians \cite{GGU}. All this can be be put into
one diagram called the \emph{Tulczyjew triple}:
\begin{equation}\xymatrix@C-5pt{
 &\sT^\ast E^\ast \ar[rrr]^{\widetilde{\Pi}_\ze}
\ar[ddl]_{\pi_{E^\ast}} \ar[dr]^{\sT^\ast\pi}
 &  &  & \sT E^\ast \ar[ddl]_/-25pt/{\zt_{E^\ast}} \ar[dr]^{\sT\pi}
 &  &  & \sT^\ast E \ar[ddl]_/-25pt/{\sT^\ast\zt} \ar[dr]^{\pi_E}
\ar[lll]_{\varepsilon}
 & \\
 & & E \ar[rrr]^/-20pt/{\zr}\ar[ddl]_/-20pt/{\zt}
 & & & \sT M\ar[ddl]_/-20pt/{\zt_M}
 & & & E\ar[lll]_/+20pt/{\zr}\ar[ddl]_{\zt}
 \\
 E^\ast\ar[rrr]^/-20pt/{id} \ar[dr]^{\pi}
 & & & E^\ast\ar[dr]^{\pi}
 & & & E^\ast\ar[dr]^{\pi}\ar[lll]_/-20pt/{id}
 & & \\
 & M\ar[rrr]^{id}
 & & & M & & & M\ar[lll]_{id} &
}\label{F1.3b}\end{equation} 

\medskip\noindent The left-hand side is Hamiltonian,
the right-hand side is Lagrangian, and the phase dynamics lives in the
middle.
Note finally that the above formalisms can still be generalized to
include constraints (cf. \cite{GG1}) and that a rigorous optimal
control theory on Lie algebroids can be developed as well
\cite{CM,GJ}.
\section{Kirillov brackets and QD-algebroids}
From the geometric point of view, of a particular interest are
brackets on the spaces of sections of vector bundles given in
differential terms. As examples we can consider the Lie bracket of
vector fields (as sections of $\sT M$) and the Poisson (or
Legendre) bracket on $C^\infty(M)$ (viewed as the space of
sections of the trivial bundle $M\ti\R\to M$) for a symplectic
(resp., contact) manifold $M$.

In [Ki], Kirillov  introduced \emph{local Lie algebra} brackets on
line bundles over a manifold $M$ as Lie brackets
on their sections given by local operators. These brackets we will
call \emph{Kirillov brackets}. According to the  Peetre
Theorem \cite{JP}, local operators are locally differential
operators, so we can as well deal locally with brackets defined by
bi-differential operators.

The fundamental fact discovered in [Ki] is that these operators
have to be of the first order and then, locally, they reduce to
the conformally symplectic Poisson or Lagrange brackets on the
leaves of a certain generalized foliation of $M$. For the trivial
bundle, i.e.  for the algebra $C^\infty(M)$ of functions on $M$,
the local brackets reduce to Jacobi brackets. Hence, the line
bundles equipped with a Kirillov bracket are sometimes called
\emph{Jacobi bundles}.

\begin{theorem} Any Kirillov bracket on sections of the trivial bundle $M\ti\R$ (i.e.  on $C^\infty(M)$) is a Jacobi bracket.
\end{theorem}

Note that in the above we view $C^\infty(M)$ as a $C^\infty(M)$-module, not as an algebra! A pure algebraic version of the above result is also
valid \cite[Theorem 4.2]{Gr}. In the purely algebraic context, we
replace the algebra $C^\infty(M)$ with an associative commutative
algebra $\A$, and the space $\Sec(E)$ of sections of a vector
bundle $\zt:E\to M$ with an $\A$-module $\cE$. We can define
linear differential operators $D:\cE_1\to\cE_2$ between two
$\A$-modules as follows. Let us observe first that, for
$f\in\A$, we can construct a new operator $\zd(f)D:\cE_1\to\cE_2$
as the `commutator' $[D,m_f]$:
\be (\zd(f)D)(x)=D(fx)-fD(x).
\end{equation}

\begin{definition} We say that a $\k$-linear operator
$D:\cE_1\to\cE_2$ between two $\A$-modules is a \emph{differential
operator of order} $\le n$ if
\begin{equation}\label{dop}
\zd(f_0)\cdots\zd(f_{n})D=0 \quad \text{for all}\quad
f_0,\dots,f_{n}\in\A\,.
\end{equation} The set of all such linear differential operators will be denoted $\Diff_n(\cE_1;\cE_2)$, or simply
$\Diff_n(\cE)$ if $\cE_1=\cE_2=\cE$.
\end{definition}

\noindent Note that if $\k$ is of characteristic 0, the condition
(\ref{dop}) can be replaced by
\be\zd(f)^{n+1}D=0 \quad \text{for
all}\quad  f\in\A\ee
and that the idea of defining differential operators in this pure algebraic context goes back to Grothendieck and Vinogradov \cite{Vi}.
It is easy to see that a zero-order
differential operator $D$ is just a module homomorphism, i.e.  an
$\A$-linear map.

\medskip\noindent
{\bf Problem.} Prove that, for an associative commutative algebra
$\A$ with unit $\1$, any first-order differential operator $D:\A\to\A$ is of the form
\be D(g)=X(g)+fg\ee
for a certain $f\in\A$ and $X\in\Der(\A)$.

\medskip
In $\Diff_0(\cE)$ there is a special class, $\A_\cE$, of
zero-order differential operators which are just
multiplications $m_f$ by elements $f$ of $\A$. Hence, in $\Diff_1(\cE)$
there is a special class, $\QD(\cE)$, of operators $D$ such that
$\zd(f)D\in\A_\cE$ for all $f\in\A$. We call them \emph{derivative endomorphisms},
\emph{quasi derivations}, or \emph{covariant differential
operators}. In other words, a $\k$-linear operator $D:\cE\to\cE$ is
a derivative endomorphism if and only if, for all $f\in\A$, there
is $\wh{D}(f)\in\A$ such that
\begin{equation}\label{qd}
D(fX)=fD(X)+\wh{D}(f)X
\end{equation}
for all $X\in\cE$. Of course, if $\cE=\A$ is the trivial module, any quasi derivation is actually a first-order differential operator on the algebra $\A$.

\medskip\noindent
{\bf Problem.} Let $D,D_1,D_2\in\QD(\cE)$. Prove that the
commutator $[D_1,D_2]$ is again in $\QD(\cE)$, that the map
\be\wh{D}:\A\to\A\quad f\mapsto\wh{D}(f)\ee is a derivation, and
that $D\mapsto \wh{D}$ is a homomorphism of the Lie algebra
$\QD(\cE)$ with the commutator bracket into the Lie algebra
$\Der(\A)$ of derivations of $\A$. We call this map the
\emph{universal anchor map}.

\medskip
For multilinear operators we define analogously the corresponding commutators with respect to the $i$'th variable,
\be
\zd_{i}(f)D(x_1,\dots,x_p)=D(x_1,\dots,fx_i,\dots,x_p)-fD(x_1,\dots,x_p)\,,
\end{equation}
and call a multilinear operator $D$ to be \emph{of order}
$\le n$ if $\zd_{i}(f_0)\cdots\zd_{i}(f_n)D=0$ for all $f_0,\dots
f_n\in\A$ and all $i$. This actually means that $D$ is of order
$\le n$ with respect to each variable separately.

\begin{definition} \emph{A differential Loday bracket} on $\A$
is a Loday  bracket on $\A$ given by a bi-differential
operator.
\end{definition}

\begin{proposition} [\cite{GMbin}] If $\A$ has no nontrivial nilpotent elements,
then every differential Loday bracket on $\A$ is actually of the order $\le 1$,
thus a Jacobi bracket.
\end{proposition}

\subsection{QD-algebroids}
A non-trivial differential requirement for a bracket
$[\cdot,\cdot]$ on sections of a vector bundle $\zt:E\to M$ is
that the bracket is a quasi-derivation with respect to each
variable separately, i.e.  $[X,\cdot]$ and $[\cdot,X]$ are
quasi-derivations for each $X\in\cE=\Sec(E)$. Hence,
$[X,fY]=f[X,Y]+\zr(X)(f)Y$    and $[fX,Y]=f[X,Y]-\zs({Y})(f)X$ for
all  $X,Y\in\cE$  and  all $f\in\A=C^\infty(M)$, where
$X\mapsto\zr(X)\in\sT M$ and $Y\mapsto\zs(Y)\in\sT M$ is,
respectively, the \emph{left} and the \emph{right anchor map}. Such brackets we
will call \emph{QD-algebroid brackets}. All skew algebroid
brackets and all Poisson brackets are of this type. The difference
is that the anchor map is of order 0 for skew algebroids and of
order 1 for a Poisson structure ($f\mapsto\zr(f)$ is just passing
to the Hamiltonian vector field).
If  the  anchor maps are of the order 0 ($\A$-linear),
we  speak  about an \emph{algebroid} (cf. \cite{GU2}).
One can prove the following, somehow unexpected, fact.

\begin{theorem} [\cite{He,QD}] \ Every QD-algebroid of rank $>1$ is an algebroid. In other words, the anchor maps may be of the order 1 on line bundles only.
\end{theorem}

\par
If the bracket is a Lie bracket, we will speak about  a
\emph{Lie QD-algebroid} (resp., \emph{Lie algebroid}). We can also
consider \emph{Leibniz QD-algebroid} (resp., \emph{ Leibniz
algebroid}) requiring additionally only the Jacobi identity
(\ref{J}) without the skew-symmetry assumption. We will not use
the term \emph{Loday algebroid} in this case to avoid a confusion
with another concept of a Loday algebroid \cite{GKP}.

Actually, there is no big difference between Leibniz and Lie
QD-algebroids, as the Jacobi identity forces the skew-symmetry. In
particular, the left anchor must be equal to the right anchor. We can sum up these results as follows.

\begin{theorem} [\cite{QD,GMbin}]\label{tx} Let $E$ be a vector bundle over $M$ and let $[\cdot,\cdot]$
be a bilinear bracket operation  on  the  $C^\infty(M)$-module
$\cE=Sec(E)$ which satisfies the Jacobi identity (\ref{J}) and which
is a quasi-derivation with respect to both arguments.
\begin{description}
\item{(a)} If $rank(E)>1$, then there is a vector bundle  morphism
$\zr:E\ra TM$ over the identity map on  $M$  such  that
$\zr([X,Y])=[\zr(X),\zr(Y)]$ and
\be
[fX,gY]=fg[X,Y]+f\zr(X)(g)Y-g\zr(Y)(f)X\,,
\ee
for all $X,Y\in\cE$,
$f,g\in C^\infty(M)$. Moreover, $[X,Y](p)=-[Y,X](p)$ if $\zr_p\ne
0$. \item{(b)} If $rank(E)=1$, then the bracket is skew-symmetric
and defines a Kirillov bracket which,  locally,  is  equivalent  to a
Jacobi bracket (\ref{Jabra}).
\end{description}
\end{theorem}

\begin{corollary}
Lie QD-algebroids on $E$ are exactly Lie algebroids if
$rank(E)>1$, and local  Lie  algebras  in  the  sense of Kirillov
if  $rank(E)=1$.
\end{corollary}

\begin{corollary}
A Lie algebroid on a vector bundle  $E$ of rank $>1$ is just a Lie
bracket on sections of $E$ which is a quasi-derivation with
respect to one (hence both) argument.
\end{corollary}

\section{Graded manifolds}


\begin{definition}
By a \emph{graded manifold} we will understand a supermanifold $\cM$
with a $\N^k$-gradation in the structure sheaf that agrees with
the parity. This means that, for any weight $w\in\N^k$, the
homogeneous functions of weight $w$ have parity coinciding with
the parity of the total weight $w=w_1+\cdots w_k$. For $k=1$,
one can also think that there is an atlas whose local coordinates
have integer weights: odd coordinates have odd weights, and even
coordinates even weights, that are preserved by changes of
coordinates. An \emph{$\N^k$-manifold of degree $d\in\N$} is a
$\N^k$-graded manifold whose local coordinates have total weights
$\le d$. A \emph{symplectic manifold of degree $r\in\N^k$} is an
$\N^k$-graded manifold equipped with a homogeneous symplectic form
of degree $r$.
\end{definition}

\begin{remark} There are various, also more general, concepts of a graded manifold, but the above will be sufficient for our purposes.
We will assume in this note that the graded manifolds are
\emph{complete}, i.e.  the even coordinates of non-zero weights
take all real values. This is to avoid considering, for instance,
open subsets in vector spaces instead of the whole vector spaces.
Basic concepts and facts concerning $\Z$-graded manifolds can be found
in \cite{Vo}. Note that $\N$-manifolds (called also
\emph{N-manifolds}) have been first studied by \v Severa
\cite{Sev} and Roytenberg \cite{Roy}.
\end{remark}

Let us remark that the $\N^k$-grading can be conveniently
encoded by means of the collection of {\em weight vector fields}
which are jointly {\it diagonalizable}, i.e.  there is an atlas of
charts with local coordinates $(x^a)$ in which
\be\label{euler}
\zD^s=\sum_aw^s_ax^a\pa_{x^a}\,,\quad s=1,\dots,n\,,
\ee
where $w^s_a=w^s(x^a)\in\N$. An $\N^k$-manifold is complete if and only
if each weight vector field is complete, i.e.  induces an action of the multiplicative $\R$.

\begin{example} If $\zt:E\ra M$ is a vector bundle, then $E[d]$ is an N-manifold of degree $d$, if we consider the
basic functions being of weight $0$ and functions linear in fibers
being of degree $d$. That the coordinate changes preserve the
weights is equivalent to preserving the vector bundle structure.
Thus, every $\N$-manifold of degree $1$ is of the form $E[1]$ with
the algebra of smooth functions $C^\infty(E[1])=\Gr(E^\*)$. The
corresponding weight vector field is the Euler vector field.
\end{example}

\begin{remark} As the Grassmann algebra $\Gr(E^*)$ can be understood as the algebra of smooth functions on the graded manifold $E[1]$ (an {\it N-manifold of degree} 1 in the terminology of \v Severa and Roytenberg \cite{Roy, Sev}),
following \cite{Va,Va1} we can view the de Rham derivative $\xdp$ as a vector field of degree 1 on $E[1]$. This vector field is {\it homological}, $(\xdp)^2=0$, if and only if we are actually dealing with a Lie algebroid.
In local supercoordinates $(x,\y)$ associated canonically with our standard affine coordinates $(x,y)$, we have
\be\label{sdR} \xdp=\frac{1}{2}c^k_{ij}(x)\y^j\y^i\pa_{\y^k} + \zr^b_i(x)\y^i\partial _{x^b}\,.
\ee
\end{remark}

\begin{example}({\bf Symplectic $\N$-manifolds of degree 1})
A {symplectic manifold of degree 1} is an $\N$-manifold of degree
$1$, thus $E[1]$ for a vector bundle $E$ over $M$, equipped with a
symplectic form of degree $1$. It is easy to see that in this case
$E$ has to be linearly symplectomorphic to the cotangent bundle $\sT^\ast
M$ with the canonical symplectic form.  The corresponding
(super)Poisson bracket $\{\cdot,\cdot\}$ on $C^\infty(T^\ast[1]
M)$ coincides with the Schouten bracket of multivector fields on
$M$. Functions of degree $2$ (quadratic) correspond to bivector
fields $\Lambda$ on $M$. Moreover, the bracket is odd and the
homological condition, $\{\Lambda,\Lambda\}=0$, means that $\Lambda$
is a Poisson tensor. The {derived bracket}, $\{ f,g\}_\Lambda=\{
\{ f,\Lambda\},g\}$, is closed on basic functions where it
coincides (up to a sign) with the Poisson bracket of $\Lambda$.
\end{example}



\subsection{The Big Bracket}  Let $\zt:E\ra M$ be a vector bundle and let $E[1]$ be the
corresponding $\N$-manifold of degree $1$.  We will use local
coordinates $(x^a,y^i)$ in $E[1]$, where $(x^a)$ are local
coordinates (of weight $0$) in a neighbourhood $W\subset M$ and
$y^i$ are linear functions in $\zt^{-1}(W)\subset E$ (of weight
$1$), corresponding to a basis of local sections of the dual
bundle $E^\ast$.
As we know, $\sT^\ast E$ is canonically a {double vector bundle}
isomorphic to $\sT^\ast E^\ast$, with the second
bundle structure being $\sT^\ast E\ra E^\ast$. In consequence,
$\sT^\ast E$ is canonically $\N^2$-graded, with local coordinates
$(x^a,y^i,p_b,\zx_j)$ having bi-degrees $(0,0),(1,0),(1,1)$, and
$(0,1)$, respectively. This $\N^2$-grading comes from the canonical
$\Z^2$ grading, in which $(x^a,y^i,p_b,\zx_j)$ have bi-degrees
$(0,0),(1,0),(0,0), (-1,0)$, from the degree shift
by $(1,1)$ in the fibers of $\pi_E:\sT^\* E\to E$. The variables $(x,
p)$ are even and the variables $(y,\zx)$ are odd. The canonical
symplectic form has the bi-degree $(1,1)$. The corresponding
graded symplectic manifold $\cM$ we will denote $\sT^\*[(1,1)]E[1]$. The
double vector bundle structure yields canonical projections
$\zt_0:\sT^\ast E\ra E\oplus_M E^\ast$ and
\begin{equation}\label{projection}\overline{\zt}_0:\sT^\ast[(1,1)]
E[1]\simeq \sT^\ast[(1,1)] E^\ast[1]\ra (E\oplus_M E^\ast)[1]\,.
\end{equation}
Of course, any double vector bundle is also an $\N$-graded
manifold of degree 2. The corresponding weight vector field is the
sum of the two commuting Euler vector fields. Therefore, the cotangent
bundle $\sT^\ast[2] E[1]\simeq \sT^\ast[2] E^\ast[1]$ is
canonically an $\N$-manifold of degree 2 with local coordinates
$(x^a,y^i,p_b,\zx_j)$ of degrees $0,1,2,1$, respectively. In
particular, $x^a,p_b$ are even coordinates and $y^i,\zx_j$ are odd
coordinates. The canonical symplectic form \be\zw=\xd p_a\xd
x^a+\xd\zx_i\xd y^i=-\xd x^a\xd p_a+\xd y^i\xd \zx_i\ee is
homogeneous of degree 2. In both above cases, the corresponding
graded Poisson bracket of degree $(-1,-1)$ (resp., $-2$), called
sometimes the \emph{big bracket}, is completely characterized
locally by
\bea\nn &\{ p_b,x^a\}=-\{ x^a,p_b\}=\zd_b^a\,,\quad &\{\zx_j,y^i\}=\{ y^i,\zx_j\}=\zd_j^i\,,\\
&\{ p_b,y^i\}=\{ p_b,\zx_j\}=0\,,\quad &\{ x^a,y^i\}=\{
x^a,\zx_j\}=0\,.
\eea
Note that, for a vector space $V$, the big
bracket on $\sT^\*[2]V[1]=(V\oplus V^\*)[1]$ has been considered
already by Kostant and Sternberg \cite{KSt}.

The projection (\ref{projection}) induces embeddings of the
algebras of smooth functions $C^\infty(E[1])=\Gr(E^\ast)$ and
$C^\infty(E^\ast[1])=\Gr(E)$ into $ C^\infty(\sT^\*[(1,1)]E[1])$
as functions of bi-degrees $(\bullet,0)$ and $(0,\bullet)$
respectively. Moreover, functions of the total degree $1$ on $\cM$
correspond to sections of $E\oplus_ME^\ast$.

Also the $\Gr(E^\*)$-module $\zF_1(E)=\Gr(E^\*)\ot_{C^\infty(M)}\Sec(E)=\Sec(\we
E^\*;E)$ is therefore interpreted as spanned by functions of
bi-degrees $(n,1)$ on $\cM=\sT^\ast[(1,1)]
E[1]$, with $n\ge 0$. In local coordinates, elements of
$\Sec(\we^kE^\*)$ are represented by polynomials
\be\sum_{i_1<\cdots<i_k}f_{i_1\cdots i_k}(x)y^{i_1}\cdots
y^{i_k}\,, \ee and elements of $\zF_1^k(E)$ by polynomials
\be\sum_{j,i_1<\cdots<i_k}g_{i_1\cdots i_k}^j(x)y^{i_1}\cdots
y^{i_k}\zx_j\,.\ee As the canonical symplectic bracket is of the
bi-degree $(-1,-1)$, it is closed on $\zF_1(E)$ which is therefore
a (graded) Lie subalgebra of $C^\infty(\cM)$.

\medskip\noindent
{\bf Problem.} Check that the big bracket restricted to $\zF_1(E)$
is exactly the (generalized) Nijenhuis-Richardson bracket.

\begin{remark}
The big bracket is also closed on
$\Gr(E)\ot_{c^\infty(M)}\Gr(E^\*)=\Gr(E\oplus_M E^\*)$. It
coincides (cf. \cite{KS}) with a bracket considered by Buttin
\cite{But}. She considered the commutator bracket of graded
differential operators $i_K$ on $\Gr(E^\*)$ associated with
elements $K$ of $\Gr(E\oplus_M E^\*)$ by \be i_{\zm\ot
X}(\zn)=\zm\we i_X\zn\,.\ee
\end{remark}

\subsection{The de Rham derivative as a homological vector field}
Since, for $E=\sT M$, the algebra of smooth functions
$C^\infty(\sT[1]M)$ is the algebra $\zW(M)$ of
differential forms, the de Rham derivative $\xd$, being a
derivation in $\zW(M)$, represents a vector field on $\sT[1]M$. In local
coordinates $(x^a,\dot x^b)$ in $\sT[1]M$ (here, $x^a$ are even
and $\dot x^b$ are odd), \be\xd=\dot x^a\pa_{x^a}\,.\ee This
vector field is odd, so that $[\xd,\xd]_{\sT[1]M}=2\xd^2$, and
{homological}, i.e.  $[\xd,\xd]_{\sT[1]M}=0$. Its lift to the
cotangent bundle $\sT^\ast[(1,1)]\sT[1]M$, with local coordinates
$(x,\dot x,p,\zp)$ of bi-degrees $(0,0),(1,0),(1,1)$, and $(0,1)$ (and the total weights $(0,1,2,1)$)), respectively,  reads
\be\H({\xd})=p_a\pa_{\zp_a}+\dot x^a\pa_{x^a}\,.
\ee
This is a
Hamiltonian vector field with the cubic Hamiltonian of the
bi-degree $(2,1)$,
\begin{equation}\label{Hamcan}H_{\xd}=\dot x^ap_a\,.
\end{equation}
More generally, if $E$ is a Lie algebroid associated with a linear
Poisson tensor
\be\Pi =\frac{1}{2}c^k_{ij}(x)\zx_k
\partial _{\zx_i}\we \partial _{\zx_j} + \zr^b_i(x) \partial _{\zx_i}
\wedge \partial _{x^b}\,,
\ee
then we can view the algebroid de Rham
derivative $\xdp$ as a vector field of degree 1 on $E[1]$. This
vector field is \emph{homological}, $(\xdp)^2=0$ \cite{Va,Va1}. In local
supercoordinates $(x,y)$ associated canonically with our standard
affine coordinates, we have
\begin{equation}\label{sdR1} \xdp=\frac{1}{2}c^k_{ij}(x)y^jy^i\pa_{y^k} + \zr^b_i(x)y^i\partial _{x^b}\,.
\end{equation} The corresponding Hamiltonian of the lifted vector field reads
\begin{equation}\label{sHam}
H=H_\xdp=\frac{1}{2}c^k_{ij}(x)y^jy^i\zx_k + \zr^b_i(x)y^ip_b\,.
\end{equation} The Hamiltonian is homological
\be\{ H,H\}=0\,,\ee and $\xdp$ is of degree $(1,0)$, so it defines
the corresponding cohomology which can be restricted to any complex
$\A^{(\bullet,n)}$. On $\A^{(\bullet,0)}=\Gr(E^\*)$, this
cohomology is the classical Lie algebroid cohomology. On the
other hand, for any section $X$ of $E$, interpreted as an element
in $\A^{(0,1)}$, the function $\{ H,X\}$ is of degree $(1,1)$ and
represents a linear vector field $\dd_\sT^\Pi X$ on $E$, the
\emph{complete lift} of $X$  (cf. (\ref{eqn:tan_lift})) or, with a different interpretation, a linear vector field $\H^\Pi(X)$
on $E^*$, the \emph{dual complete lift} of $X$.

\subsection{The Fr\"olicher-Nijenhuis bracket revisited}
Theorem \ref{tsb} implies immediately the following.

\begin{proposition} The bracket derived from the big bracket and the Hamiltonian $H$ of $\xd^\zP$,
\be\label{debr}\{ K,L\}^{H}=\{\{ K,H\},L\}\,,\ee is closed on $\Gr(E)$ and
coincides there with the generalized Schouten bracket
$\lna\cdot,\cdot\rna^\zP$.
\end{proposition}

It is easy to see that the derived bracket (\ref{debr}) is closed also on vector-valued forms, i.e.  on
$\zF_1(E^\*)$. However, it gives not the Fr\"olicher-Nijenhuis
bracket, since it is not skew-symmetric. The Fr\"olicher-Nijenhuis
bracket differs from the derived one by a coboundary term (cf. (\ref{Fro2})).

A tensor $N\in\zF_1^1(E)$ we will call an (algebroid) \emph{pseudo-Nijenhuis tensor} if $\{ N,N\}^{H}=0$.

\begin{theorem}[cf. \cite{CGMc, Gra3}]
For $K\in\zF_1^k(E^\*)$ and $L\in\zF_1^l(E^\*)$, we have
\begin{equation}\label{fnij} [K,L]^{FN}=\{ K,L\}^{H}+(-1)^{k(l+1)}\{
i_LK,H\}\,.
\end{equation}
Any pseudo-Nijenhuis tensor $N$ is weak-Nijenhuis, $\{ H,N\}$ is a homological Hamiltonian, and
the contracted bracket (\ref{Nij cont}), corresponding to $\{ H,N\}$,
is again a Lie algebroid bracket.
\end{theorem}

\noindent For a discussion of brackets associated with the big bracket we refer to the
survey article \cite{KS}.

\section{Courant bracket and Dirac structures} Recall that if
$(M,\omega)$ is an $2n$-dimensional symplectic manifold, then the
symplectic form $\omega=\frac{1}{2}\zw_{ij}\xd x^i\we \xd x^j$
induces a vector bundle isomorphism \be\tilde{\omega}:\sT M\ni
V\mapsto -i_V\omega\in\sT^\ast M\,.\ee The inverse map
\be\tilde{\zL}=\wt{\omega}^{-1}:\sT^\ast M\ra \sT M\ee
corresponds to a Poisson tensor
$\zL=\frac{1}{2}\zL^{ij}\pa_{x^i}\we\pa_{x^j}$ \emph{via} $\wt{\zL}(\za)=i_\za\zL$.
The fact that $\zw$ is closed, $\xd\zw=0$, reads in coordinates as
\begin{equation}\label{closed}(\xd\zw)_{kij}=\frac{\pa\zw_{ij}}{\pa x^k}+
\frac{\pa\zw_{jk}}{\pa x^i}+\frac{\pa\zw_{ki}}{\pa x^j}=0\,,\end{equation}  or, equivalently,
\begin{equation}\label{Poisson}[\zL,\zL]^{kij}=\frac{\pa\zL^{ij}}{\pa x^l}\zL^{lk}+
\frac{\pa\zL^{jk}}{\pa x^l}\zL^{li}+\frac{\pa\zL^{ki}}{\pa x^l}\zL^{lj}=0\,.
\end{equation}
Note that both equations, $\xd\zw=0$ and $[\zL,\zL]=0$, equivalent
for invertible $\zL=\zw^{-1}$, make sense for an arbitrary 2-form
$\zw$ and any bivector field $\zL$ separately.

\medskip
The (common) graph of $\zw$ and $\zL=\zw^{-1}$ is a vector
subbundle $L$ of the \emph{Pontryagin bundle} $\cT M=\sT M\oplus_M\sT^\ast M$,
 \begin{eqnarray}\nn L_p&=&\{ (V_p+\zz_p)\in\sT_p M\oplus\sT^\ast_p M:V_p=\wt{\zw}(\zz_p)\}\\
 &=&\{ (V_p+\zz_p)\in\sT_p M\oplus\sT^\ast_p M:\zz_p=\wt{\zL}(V_p)\}\,.\label{D}
 \end{eqnarray}
The skew-symmetry of $\zw$ (or $\zL$) means that $L$ is isotropic
with respect to the canonical symmetric pseudo-Riemannian metric $\bk{\cdot}{\cdot}_+$ on $\cT M$, where
 \begin{equation}\label{metric} \bk{V_p+\zz_p}{U_p+\zh_p}_\pm=
 \frac{1}{2}\left(\pr{\zz_p}{U_p}\pm\pr{V_p}{\zh_p}\right)\,.\end{equation}
 The condition $\xd\zw=0$ (or $[\zL,\zL]=0$) means that $L$ is
 involutive with respect to the \emph{Courant bracket} \cite{Co1} on $\cT M$ defined by
 \begin{equation}\label{Cou} [V+\zz,U+\zh]_C=[V,U]+\left(\Ll_V\zh-\Ll_{U}\zz+\xd\bk{V+\zz}{U+\zh}_-\right)\,.
 \end{equation}
An important observation is that on any isotropic subbundle $L$
the Courant bracket coincides with the \emph{Dorfman bracket}
\cite{Do} given by
\begin{equation}\label{CD}[V+\zz,U+\zh]_D=[V,U]+\left(\Ll_V\zh-i_{U}\xd\zz\right)\,.
\end{equation}
Starting with a bivector field $\zL$ and denoting, for arbitrary
functions $f,g$, the corresponding Hamiltonian vector fields
$V_f,V_g$, respectively, we get \be[V_f+\xd f,V_g+\xd
g]_D=[V_f,V_g]+\xd\{ f,g\}\,,\ee
 so that involutivity means $V_{\{ f,g\}}=[V_f,V_g]$,
that is equivalent to $[\zL,\zL]=0$.

The Dorfman bracket is not skew-symmetric but it satisfies the
{Jacobi identity}, so it is a Loday bracket. The Courant bracket
is skew-symmetric but it does not satisfy the Jacobi identity;
there is a \emph{Jacobi anomaly}.  Both brackets coincide on any
isotropic subbundle $L$ and give a Lie algebroid bracket on the space $\cL$
of its sections if $L$ is involutive. Actually,
the Dorfman bracket is a derived bracket.  Namely, we use the
Hamiltonian (\ref{Hamcan}) to define a derived bracket out of the
canonical Poisson bracket on $\sT^\ast[2]\sT[1]M$:
\be\lna
A,B\rna=\{ A,B\}^{H_\xd}=\{\{ A,H_\xd\},B\}\,. \ee
This bracket is
of degree $-1$, so it is closed on functions of degree $1$, thus
sections of $\cT M=\sT M\oplus_M\sT^\ast M$, where it coincides
with the Dorfman bracket.

\begin{definition} A \emph{Dirac structure} is a maximal isotropic and
involutive subbundle $L$ of $(\cT M, [\cdot,\cdot]_D)$.
We call a vector field $V$ on $M$ an {$L$-Hamil\-tonian vector field} with an {$L$-Hamiltonian function} $f$ if $(V+\xd f)\in L$.
\end{definition}

Let $\fH$ be the set of all $L$-Hamiltonian vector fields, $\fH_0$
be the set of all $L$-Hamiltonian vector fields with the
Hamiltonian $0$, and $\fA$ be the set of projections of sections
of $L$ onto $\sT M$.

\begin{theorem}
The families $\fH$, $\fH_0$, and $\fA$ are Lie algebras of vector
fields. Moreover, there is a canonical Poisson bracket on the
space $\fH$ of all Hamiltonians,
\be\{ f_1,f_2\}_\cL=\pr{V_1}{\xd
f_2}\quad \text{if} \quad (V_i+\xd f_i)\in\cL\,,
\ee
that endows this space with a
Lie algebra structure. If \ $V_i$ is an $L$-Hamiltonian vector field with
an $L$-Hamiltonian $f_i$, $i=1,2$, then $\{ f_1,f_2\}_\cL$ is an
$L$-Hamiltonian of the $L$-Hamiltonian vector field $[V_1,V_2]$.
\end{theorem}

\noindent Dirac structures induce presymplectic foliations on $M$ as
follows.

\begin{theorem}
The Lie algebra of vector fields $\fA$ induces a (generalized)
foliation $\cF$ of $M$. Every leaf $Y$ of this foliation is a
presymplectic manifold with the closed two-form $\zw_Y$ induced
from the map
\be\fA\ti\fA\ni(V_1,V_2)\mapsto\zW_L(V_1,V_2)=\zz_1(V_2)\,.\ee
Here, $\zz_1$ is any 1-form satisfying $(V_1+\zz_1)\in\cL$.
Moreover, $L$-Hamiltonians are functions constant along the
characteristic distributions of these presymplectic forms and the
corresponding $L$-Hamiltonian vector fields are their Hamiltonian
vector fields with respect to the presymplectic forms.
\end{theorem}

\begin{example} {\bf(Dirac constrains)}
Dirac structures on manifolds provide a geometric setting for
Dirac's theory of constrained mechanical systems.
Let $Y\subset M$ be a submanifold determined by $r$ independent
constraints
\be\phi_1(p)=\cdots=\phi_r(p)=0\,.\ee
Note that the map $\phi=(\phi_i)$ defines actually a
foliation $\cF=\{\phi=const\}$, not a single submanifold.
Let $\sT\cF\subset\sT M$ be the corresponding distribution and
$(\sT\cF)^0\subset\sT^\ast M$ its annihilator (spanned by
$\xd\phi_i$). The collective constraint $\phi$ defines a Dirac
structure $L^\phi\subset\cT M$ with the fibers
\be
L_p^\phi=\{(V_p+\zz_p)\in\cT M:V_p\in\sT_p\cF\ \text{and\ } \
\zz\in\wt{\zw}(V_p)+(\sT\cF)^0_p\}\,.
\ee
In this case, $L^\phi$-Hamiltonians are functions $f$ satisfying
\be\{ f,\phi_i\}=\mu^j\{\phi_j,\phi_i\}\,,\quad \mu^j\in C^\infty(M)\,.
\ee
If $\phi$ consists of \emph{first-class constraints},
$\{\phi_i,\phi_j\}=0$, then $L^\phi$-Hamiltonians are \emph{first-class
functions}, $\{ f,\phi_i\}=0$, and the Poisson bracket on the
algebra of first-class functions is the original symplectic
Poisson bracket.
The bracket of first-class functions is again first-class:
\be\{
\{ f,g\},\phi_i\}=\{\{ f,\phi_i\},g\}+\{ f,\{ g,\phi_i\}\}=0\,.\ee
If $\phi$ are \emph{second-class
constraints}, i.e.  the matrix  $\left(\{\phi_i,\phi_j\}\right)$ is
invertible, $\left(\{\phi_i,\phi_j\}\right)^{-1}=(c^{ij})$, then
it defines a foliation into symplectic submanifolds, so any function
is $L^\phi$-Hamiltonian and
the Poisson bracket on the algebra of $L^\phi$-Hamiltonians is
the \emph{Dirac bracket}
\be\{ f,g\}_{\cL^\phi}=\{ f,g\}-\{ f,\phi_i\}
c^{ij}\{\phi_j,g\}\,.\ee
\end{example}
\subsection{Multi-Dirac and Poly-Dirac structures}

The Dorfman bracket (\ref{CD}) can be immediately generalized (cf.
\cite{BS}) to a bracket on sections of $\cT^\bullet M=\sT
M\oplus_M\wedge^\bullet\sT^\ast M$, where
\be\wedge^\bullet\sT^\ast M=\bigoplus_{k=0}^\infty\wedge^k\sT^\ast
M\,,
\ee
so that sections of $\cT^\bullet M$ are of the form $(X+\zw)$, where $X$ is a vector field and $\zw$ is a differential form.
The bracket, which we will call
the \emph{Grassmann-Dorfman bracket}, is formally given by the
same formula $(\ref{CD})$ and it is also a Loday bracket. It can
be reduced to a bracket $\lna\cdot,\cdot\rna^n$ on sections of the
\emph{Pontryagin Bundle of degree $n$}, i.e.  the bundle
$\cT^nM=\sT M\oplus_M\we^n\sT^\ast M$, $n\in\mathbb{N}$, the
\emph{Grassmann-Dorfman bracket of degree $n$}, being an example
of a \emph{Loday algebroid bracket} \cite{GKP}. In particular, the
projection $\zr:\cT M\to\sT M$ onto the first summand yields
the left anchor of the bracket.

Note that the Grassmann-Dorfman bracket is a part of the
\emph{graded Courant bracket} introduced in \cite{VYL} on sections
of $\we^\bullet\sT M\oplus_M\we^\bullet\sT^\ast M$. We will not
discuss the latter generalization closer, as the Grassmann-Dorfman
bracket will be sufficient for our purposes.
On $\cT^\bullet M$ we have another canonical structure, namely the
non-degenerate symmetric pairing with values in
$\we^{\bullet}\sT^\* M$,
\begin{equation}\label{pairing}\la X+\zw,Y+\zh\ran=\frac{1}{2}\left(i_X\zh+i_Y\zw\right)\,,
\end{equation}
where $\zw,\zh\in\zW(M)$. This pairing is non-degenerate also on every $\cT^nM$.

\begin{definition} A vector subbundle $L$ of the Pontryagin bundle of degree $n$
is called a \emph{multi-Dirac structure of degree $n$} if it is maximally
isotropic with respect to the above pairing and \emph{involutive}, i.e. 
whose sections are closed with respect to the Grassmann-Dorfman bracket.
\end{definition}

\noindent The following is well known (see e.g. \cite{BS}).

\begin{proposition}
The graph \be\mathfrak{G}(\za)=\{ X+i_X\za: X\in\sT
M\}\subset\cT^nM\ee of an $(n+1)$-form $\za$ on $M$ is a maximally
isotropic subbundle in $\cT^nM$. It is involutive if and only if
$\za$ is closed. The form is non-degenerate if and only the
projection of $\fG(\za)$ on the second summand is injective.
\end{proposition}

Note that closed non-degenerate $(n+1)$-forms are sometimes called
\emph{$n$-plectic} (\emph{multisymplectic}) \emph{structures}.
This justifies the following.

\begin{definition} A \emph{multi-Poisson structure of degree $n$}
is a multi-Dirac structure of degree $n$ which is the graph of a map
$\we^n\sT^\*M\supset D\to \sT M$ (on a vector subbundle domain $D$).
\end{definition}

We can slightly generalize the above concepts by considering, for
a real vector space $W$, the $W$-valued Grassmann-Dorfman bracket
and the pairing as follows. The $W$-valued Grassmann-Dorfman
bracket is defined on sections of $\cT_WM=\sT
M\oplus_M(\we^\bullet\sT^\* M\ot W)$ by
\begin{equation}\label{CDV}\lna X+\zw\ot a,Y+\zh\ot b\rna_{W}=[X,Y]+\Ll_X\zh\ot b-i_Y\,\xd\zw\ot a
\end{equation} and the $W$-valued pairing \be\la\cdot,\cdot\ran_W:\cT_W M\ti_M\cT_W M\to\cT_W M\ot W\ee by
\begin{equation}\label{pairingV}\la X+\zw\ot a,Y+\zh\ot b\ran_W=\frac{1}{2}\left(i_X\zh\ot b+i_Y\zw\ot a\right)\,.
\end{equation}
It is clear that \emph{$W$-valued poly-Dirac structure of degree $n$}
should be understood as maximal isotropic and involutive subbundles in $\cT^n_WM$.
If $W=\R^k$, we will speak about \emph{poly-Dirac structures}. An example is given
by the graph of a \emph{$W$-valued polysymplectic form} $\za\in\zW^2(M)\ot W$
(called just \emph{polysymplectic} if $W=\R^k$, cf. \cite{FG,Gu1}) which is a $W$-valued poly-Dirac structure (of degree 1).
This justifies the following definition which agrees with the concept of a
poly-Poisson structure studied in \cite{IMV}.

\begin{definition} A \emph{$W$-valued poly-Poisson structure of degree $n$}
is a $W$-valued poly-Dirac structure of degree $n$ which is the graph of a map
$\we^n\sT^\*M\ot W\supset D\to \sT M$ (on a vector subbundle domain $D$).
\end{definition}

\noindent Actually, we can replace $\sT M$ with an arbitrary Lie algebroid
$E$ and replace $\cT_W^nM$ with the Lie algebroid $W$-valued
Pontryagin bundle of degree $n$,
\begin{equation}\label{VP}\cP^n_WE=E\oplus_M\left(\we^nE^\*\ot
W\right)\,.
\end{equation}

\subsection{Courant algebroids}
Algebraic properties of the Courant bracket led to the concept of
\emph{Courant algebroid}. The original idea of Liu,
Weinstein, and Xu  \cite{LWX} was based on the observation
that $\cT M$, endowed with the Courant bracket, plays the role of a `double'
object in the sense of Drinfeld \cite{Dr} for a pair of Lie algebroids. Let us
recall that, in complete analogy with Drinfeld's Lie bialgebras, in the category of Lie algebroids there also exist `bi-objects', Lie bialgebroids, introduced by Mackenzie and Xu \cite{MX}. On the other hand, every Lie bialgebra has a double which is a Lie algebra. This is not so for
general Lie bialgebroids. Instead, Liu, Weinstein, and Xu showed that the double of a Lie
bialgebroid is a more complicated structure they call a {\it Courant algebroid}, $\cT M$ with the Courant bracket being a special case. In the general case:
\begin{itemize}
 \item the Pontryagin bundle $\cT M$ with the
canonical symmetric pairing is replaced with a vector bundle $E\ra
M$ equipped with a nondegenerate  symmetric bilinear form $( \cdot
, \cdot )$ on the bundle;

\item the Courant bracket is replaced with a  skew-symmetric
bracket $[\cdot , \cdot ]$ on $\Sec (E)$;

 \item the canonical projection $\cT M\ra\sT M$ is
replaced by a bundle map $\rho :E\ra TM$.
It induces a map $\cD :  C^{\infty}(M)\ra \Sec (E)$
defined  by $\cD = \half \beta^{-1}\rho^{*} d$, where $\beta $ is
the isomorphism between $E$ and $E^\ast$ given by the bilinear
form.  In other words, \be(\cD f, e)= \half  \rho (e) f\,. \ee
\end{itemize}

\begin{definition} [cf. \cite{LWX}] A \emph{Courant algebroid} is a vector
bundle $E\ra M$ equipped with a nondegenerate  symmetric bilinear
form
 $( \cdot , \cdot )$ on the bundle,  a  skew-symmetric
bracket $[\cdot , \cdot ]$ on $\Sec (E)$, and a bundle map $\rho
:E\ra \sT M$ (the \emph{anchor}) such that:
 \begin{enumerate}
\item For any $e_{1}, e_{2}, e_{3}\in \Sec (E)$, $[[e_{1}, e_{2}],
e_{3}]+(cyclic)=\cD T(e_{1}, e_{2}, e_{3})$, where $T(e_{1}, e_{2},
e_{3})$ is the function on the base $M$ defined by
\be T(e_{1}, e_{2}, e_{3})=\third ([e_{1}, e_{2} ], e_{3})+(cyclic)\,;
\ee
\item   for any $e_{1}, e_{2} \in \Sec (E)$, $\rho ([e_{1},
e_{2}])=[\rho (e_{1}), \rho  (e_{2})];$
\item for any $e_{1}, e_{2} \in
\Sec (E)$ and $f\in C^{\infty} (M)$, \be[e_{1}, fe_{2}]=f[e_{1},
e_{2}]+(\rho (e_{1})f)e_{2}- (e_{1}, e_{2})\cD f \,;\ee
\item
$\rho \circ \cD =0$, i.e.   for any $f, g\in C^{\infty}(M)$,
$(\cD f,  \cD g)=0$;
\item for any $e, h_{1}, h_{2} \in \Sec (E)$,
\be\rho (e) (h_{1}, h_{2})=([e , h_{1}]+\cD (e ,h_{1}) ,
h_{2})+(h_{1}, [e , h_{2}]+\cD  (e ,h_{2}) )\,.\ee
\end{enumerate}
\end{definition}

\noindent In what follows we will give equivalent `user friendly' definitions.

\subsection{Courant algebroid \emph{via} the Dorfman bracket}
For a Courant algebroid, instead of the skew-symmetric bracket
with the anomaly in the Jacobi identity, we can consider a bracket
which, like the Dorfman bracket,  is not skew-symmetric, but
satisfies the Jacobi identity, i.e.  which is a Loday bracket.
This new operation on sections of $E$ is defined by
 \be e_1\circ e_2=[e_1,e_2]+\half\cD(e_1,e_2)\,,\ee
so that the Courant bracket is the skew-symmetrization of "$\circ$",
\be[e_1,e_2]=\half\left(e_1\circ e_2-e_2\circ e_1\right)\,.
\ee
The Jacobi anomaly vanishes
and we can state an equivalent simplified definition of Courant
algebroid as follows (see \cite{GMgr,Uch}).

\begin{definition} A \emph{Courant algebroid} is a vector bundle $E\ra M$ equipped with a nondegenerate symmetric
bilinear form $({\cdot},{\cdot})$ on $E$ and a Leibniz product
(bracket) $\circ$ on $\Sec(E)$, together with a vector bundle map
(the anchor) $\zr:E\ra \sT M$, which are {compatible} with
$({\cdot},{\cdot})$, that is, \be\zr(X)( Y,Z)=( X,Y\circ Z+Z\circ Y) \ee
 and
\be\zr(X)(Y,Z)=( X\circ Y,Z)+( Y,X\circ Z)\,. \ee
\end{definition}

\noindent The latter invariance of the pairing $({\cdot},{\cdot})$ with
respect to the left multiplication implies the standard property
of the anchor map $\zr$:
\bea\nn X\circ(fY)&=&f(X\circ Y)+\zr(X)(f)Y\,.\\
\zr(X\circ Y)&=&[\zr(X),\zr(Y)]\,.
\eea
Besides its simplicity, this definition allows for considering as well the symmetric form
$({\cdot},{\cdot})$ being degenerate.
Note finally that one can define Nijenhuis tensors also for Courant algebroids \cite{CGM,Gra3,KS2} that leads to the concept of `generalized geometries' in the spirit of the Hitchin's \emph{generalized complex geometry} \cite{Hi} (see also \cite{Gu}).


\subsection{Symplectic $\N$-manifolds of degree 2}
The following characterizations of symplectic $\N$-manifolds of degree
2 and Courant algebroids as certain Hamiltonian systems are due to
Roytenberg \cite{Roy}.

\begin{theorem}
There is a one-to-one correspondence between symplectic
$\N$-manifolds of degree two, $(\cM,\zw)$, and vector bundles
$\zt:E\ra M$ equipped with a pseudo-Riemannian structure
$(\cdot,\cdot)$, i.e.  a symmetric non-degenerate two-form in
fibers. The symplectic manifold $\cM_E$ associated with
$(E,(\cdot,\cdot))$ is the pullback of $\sT^\*[2]E[1]$ with
respect to the embedding $E\hookrightarrow E\oplus_M E^\*$ given
by $X\mapsto X+(X,\cdot)$, i.e.  it completes the commutative
diagram {\[
\begin{array}{ccc}
\cM & \longrightarrow  & \sT^{*}[2]E[1]\\
\downarrow  &  & \downarrow \\
E[1] & \longrightarrow  & (E\oplus E^{*})[1]
\end{array}\]}
Moreover, the symplectic form $\zw$ is the pullback of the
canonical symplectic form on $\sT^{*}[2]E[1]$.
\end{theorem}

\begin{theorem} There is a one-to-one correspondence between Courant algebroids and symplectic
N-manifolds of degree $2$, $(\cM_E,\zw)$, equipped with a cubic
homological Hamiltonian $H$, $\{ H,H\}=0$. In this correspondence,
we identify sections of $E$ with functions of degree $1$ on
$\cM_E$, basic functions (functions on $M$) with functions of
degree $0$ on $\cM_E$, and the pseudo-riemannian metric with the
Poisson bracket, $(X,Y)=\{ X,Y\}$. The (Dorfman) algebroid bracket
on sections of $E$ is the derived bracket $X\circ Y=\{\{ X,H\},Y\}$.
\end{theorem}

Consider local coordinates $(x^a,\zz^i,p_b)$ in $\cM_E$
corresponding to coordinates $(x^a)$ on $M$ and a local basis $\{
e_i\}$ of sections of $E$ such that $(e_i,e_j)=g_{ij}=const$\,,
$e_i=g_{ij}\zz^j$ interpreted as a linear function on $E$. Then, the symplectic
form $\zw$ reads \be\zw=\xd p_a\xd x^a+\frac{1}{2}g_{ij}\xd\zz^i
\xd\zz^j\,,\ee and any cubic Hamiltonian is of the form \be
H=\zz^i\zr^a_i(x)p_a-\frac{1}{6}\zvf_{ijk}(x)\zz^i\zz^j\zz^k\,.\ee
For the corresponding Courant algebroid, the Dorfman bracket and
the anchor are uniquely determined by \be([e_i,
e_j],e_k)=\zvf_{ijk}(x)\,,\quad
\rho(e_i)=\zr^a_i(x)\partial_{x^a}\,.\ee
\section{Nambu-Poisson brackets}
There are two main ways of generalizing the notion of a Lie
algebra. One way, already discussed, is to drop the skew-symmetry assumption and consider Loday brackets.
Another concept is due to Filippov, who developed a theory of
brackets with more than two arguments, i.e.  \emph{$n$-ary
brackets}. In \cite{Fi}, he proposed a definition of such
structures which we shall call \emph{Filippov algebras}, with a
version of the Jacobi identity for $n$-arguments which we will call
\emph{Filippov identity}:
\begin{equation}\label{FI}
\{f_1,\dots,f_{n-1},\{ g_1,\dots,g_n\}\}=\sum_{k=1}^n
\{g_1,\dots,\{f_1,\dots,f_{n-1},g_k\},\dots,g_n\}.
\end{equation}
A \emph{Filippov bracket} is a skew-symmetric $n$-ary bracket satisfying (\ref{FI}). Note that in the binary case ($n=2$), the Filippov identity coincides with the Jacobi identity.
Independently, Nambu \cite{Nam}, looking for generalized
formulations of Hamiltonian Mechanics, found $n$-ary analogs of
Poisson brackets for which Takhtajan \cite{Ta} rediscovered the
Filippov identity (and called it Fundamental Identity).
This leads to the concept of a \emph{Nambu-Poisson bracket}, defined on
a commutative associative algebra, which is a Filippov bracket satisfying additionally the
\emph{Leibniz rule}: \be
\{f_1f_1',f_2,\dots,f_n\}=f_1\{f_1',\dots,\dots,f_n\}+
\{f_1,\dots,\dots,f_n\}f_1'.
\end{equation}

\begin{example} On $\R^m$,
the $n$-ary bracket operation
\begin{equation}\label{01}
\{ f_1,\dots,f_n\}={\rm det}\left(\frac{\pa f_i}{\pa x_j}\right)_{i,j\le n},
\end{equation}
where $n\le m$, is a Nambu-Poisson bracket. Actually, each nonsingular Nambu-Poisson $n$-ary bracket, with $n>2$, is
locally of this form \cite{Ga,GM,MVV}.
\end{example}

\par
It is now clear that we can combine both generalizations and
define \emph{Filippov-Loday} algebras as those which are  equipped
with  $n$-ary brackets, not skew-symmetric in general, but
satisfying the Filippov identity. We can also define a Loday
version of Nambu-Poisson algebras or rings.

\begin{definition} Let $\A$ be an associative
commutative algebra.
An {$n$-ary bracket} on $\A$ is called a \emph{Nambu-Loday} bracket if it
satisfies the Filippov identity (\ref{FI}) and the
Leibniz rule with respect to each argument $i=1,\dots,n$:
\be
\{f_1,\dots,f_if_i',\dots,f_n\}=f_i\{f_1,\dots,f_i',\dots,f_n\}+
\{f_1,\dots,f_i,\dots,f_n\}f_i'\,.
\end{equation}
\end{definition}

We encountered an unexpected phenomenon while looking for canonical examples of
Nambu-Loday brackets. One can show that, for a wide variety of associative
commutative algebras, including algebras of smooth functions, we
get nothing more than what we already know, since Nambu-Loday
brackets have to be skew-symmetric automatically. In particular, we can skip requiring the
skew-symmetry  in the standard definition of a Nambu-Poisson bracket.
Recall that we have obtained a similar negative result for a
Loday-type generalization of Lie algebroids (Theorem \ref{tx}).

\begin{theorem}[\cite{GMnas}] If $\A$ is an associative commutative algebra over a field
of characteristic 0 and $\A$ contains no nilpotents, then every
Nambu-Loday bracket on $\A$ is skew-symmetric. In particular, any
Nambu-Loday bracket on $C^\infty(M)$ is a Nambu-Poisson bracket.
\end{theorem}

\noindent For a deeper discussion of $n$-ary brackets we refer to
the review paper \cite{AI}.

\end{document}